\documentclass[12pt]{amsart}

\usepackage{amssymb}
\usepackage{latexsym}
\usepackage{verbatim}
\usepackage{fullpage}

\input {cyracc.def}
\font\ququ=cmr10 scaled \magstep1
 \font\tencyr=wncyr10 
\font\tencyi=wncyi10 
\font\tencysc=wncysc10 
\def\rus{\tencyr\cyracc}
\def\rusi{\tencyi\cyracc}
\def\rusc{\tencysc\cyracc}

\newcommand{\re}[1]{\textrm  (\ref{#1})}

\renewenvironment{proof}
{\noindent {\sl Proof.}\quad }{\hfill
$\square$ \vskip1.1ex\noindent }

\newenvironment{proof*}
{\noindent {\sl Proof.}\quad }{\hfill
$\square$}

\renewcommand{\theequation}{\thesection .\arabic{equation}}
\renewcommand{\thesubsubsection}{\theequation .\arabic{subsubsection}}

\catcode`\@=\active
\catcode`\@=11
\def\@eqnnum{\hbox to
.01pt{}\rlap{\hskip-\displaywidth(\mathbf{\theequation})}}
\catcode`\@=12

\newenvironment{s}[1]
{ \vskip1.2ex \refstepcounter{equation}
\noindent {\bf \theequation\enspace #1.} \begin{sl}}{\end{sl}
\vskip1.1ex\noindent }

\newenvironment{rem}[1]
{ \vskip1.2ex \refstepcounter{equation}
\noindent {\bf \theequation\enspace {#1}.} }{ \vskip1.1ex\noindent }

\newenvironment{subs}[1]
{\vskip1.2ex \refstepcounter{equation}
\noindent {\bf (\theequation)\quad #1.} }{\quad}

\newcommand {\be}{{\frak b}}
\newcommand {\ce}{{\frak c}}

\newcommand {\g}{{\frak g}}

\newcommand {\te}{{\frak t}}
\newcommand {\ut}{{\frak u}}

\newcommand {\gA}{{\goth A}}
\newcommand {\gC}{{\goth C}}
\newcommand {\gH}{{\goth H}}

\newcommand {\esi}{\varepsilon}
\newcommand {\ap}{\alpha}

\newcommand {\lb}{\lambda}

\newcommand {\HW}{\widehat W}
\newcommand {\HP}{\widehat\Pi}
\newcommand {\HD}{\widehat\Delta}

\newcommand {\cH}{{\mathcal H}}

\newcommand {\asst}{{\ast\,}}

\newcommand {\rk}{{\mathrm{rk\,}}}

\newcommand {\GR}[2]{{\textrm{{\bf #1}}}_{#2}}

\newcommand {\ov}{\overline}
\newcommand {\un}{\underline}

\newcommand {\Ab}{{\frak Ab}}
\newcommand {\AD}{{\frak Ad}}
\newcommand {\AN}{{\frak An}}

\newcommand {\qus}{\hfill $\square$ \vskip1.1ex\noindent}
\newcommand {\qu}{\hfill $\square$}
\newcommand {\beq}{\begin{equation}}
\newcommand {\eeq}{\end{equation}}

\newcommand{\Cat}{\mathrm{Cat}}

\newcommand{\curle}{\preccurlyeq}
\DeclareMathOperator{\rt}{rt}

\newfam\Bbbfam\newfam\eufam\newfam\eusfam%
\font\Bbbfont=msbm10 scaled 1200%
\font\olala=msam10 scaled 1200%
\font\frak=eufm10 scaled 1400%
\font\Bbbsmallfont=msbm8%
\textfont\Bbbfam=\Bbbfont\scriptfont\Bbbfam=\Bbbsmallfont
\font\euzw=eufm10 scaled 1200%
\font\euac=eufm7 scaled 1200%
\font\euacc=eufm7 scaled 1000%
\textfont\eufam=\euzw\scriptfont\eufam=\euac%
\scriptscriptfont\eufam=\euacc%
\font\euszw=eusm10 scaled 1200%
\font\eusac=eusm7 scaled 1200%
\font\eusacc=eusm7 scaled 1000%
\textfont\eusfam=\euszw\scriptfont\eusfam=\eusac%
\scriptscriptfont\eusfam=\eusacc%
\def\frak{\fam\eufam}%
\def\goth{\fam\eusfam}%
\def\Bbb{\fam\Bbbfam}%

\def\varnothing{\hbox {\Bbbfont\char'077}}
\def\square{\hbox {\olala\char"03}}

\def\cyeq{\hbox {\olala\char'064}}

\begin{document}
\setlength{\parskip}{2pt plus 4pt minus 0pt}
\hfill {\scriptsize \today} 
\vskip1ex
\vskip1ex

\title[Ideals of Heisenberg type and minimax elements]{Ideals of 
Heisenberg type and minimax elements of affine Weyl groups}
\author{\sc Dmitri I. Panyushev}
\address{Independent University of Moscow,
Bol'shoi Vlasevskii per. 11 \\
121002 Moscow, \quad Russia }
\email{panyush@mccme.ru}
\thanks{This research was supported in part by R.F.B.I. Grants no. 
01--01--00756 and  02--01--01041}
\keywords{Borel subalgebra, ad-nilpotent ideal, affine Weyl group}
\subjclass[2000]{17B20, 20F55}
\maketitle

\noindent
Let $\g$ be a simple Lie algebra 
over $\Bbb C$. 
Fix a Borel subalgebra $\be$ and a Cartan subalgebra $\te\subset\be$. 
The corresponding set of positive (resp. simple) roots is $\Delta^+$
(resp. $\Pi$). Write $\theta$ for the highest root in $\Delta^+$.
\\[.4ex] 
An ideal of $\be$ is called 
{\sf ad}-{\it nilpotent\/}, if it is contained in $[\be,\be]$.
Consequently, such a subspace is fully determined by the corresponding 
set of positive roots, and this set is called an {\it ideal\/} 
(in $\Delta^+$). 
Let ${\AD}$  denote the set of all ideals of $\Delta^+$. We regard 
$\Delta^+$ as poset with respect to the standard root order `$\curle$' 
(see Section~\ref{prelim}). Given $I\in\AD$, 
the minimal roots in it are said to be the 
{\it generators\/}. The set of generators of an ideal 
form an {\it antichain\/}
in $(\Delta^+, \curle)$, and this correspondence sets up a bijection 
between the ideals and the antichains of $\Delta^+$.
An ideal or antichain is called {\it strictly positive\/}, if it contains 
no simple roots. Another interesting class consists of Abelian ideals, i.e.,
those with the property $(I+I)\cap \Delta^+=\varnothing$.
We write $\AD_0$ (resp. $\Ab$) for the set of strictly positive 
(resp. Abelian) ideals.
Various results for $\AD$, $\AD_0$, and $\Ab$
were recently obtained in 
\cite{ath02},\,\cite{ath03},\,\cite{CP1},\,\cite{CP2},\,%
\cite{lp},\,\cite{duality},\,%
\cite{losh},\,\cite{eric}. 
\\[.4ex]
It was shown by Cellini and Papi~\cite{CP1} that there is a bijection 
between the ideals and 
certain elements of $\HW$, the affine Weyl group associated with $\g$.
Then Sommers \cite{eric} discovered a bijection between the strictly 
positive ideals (or antichains) and another class of
elements of $\HW$. Following \cite{eric}, the elements in these two classes
are said to be {\it minimal\/} and {\it maximal\/}, respectively.
Furthermore, the minimal elements of $\HW$ are in a bijection with
the points of the coroot lattice lying in a certain simplex. The same assertion 
also holds for the maximal elements (with another simplex!)
A different approach to ({\sf ad}-nilpotent) ideals relates them with the 
theory of hyperplane arrangements.
By a result of Shi \cite{shi0}, there is a 
bijection between $\AD$ and the dominant regions of the 
{\it Shi ({\rm or} Catalan) arrangement}.
It was then observed by Athanasiadis \cite{ath02} and Panyushev 
\cite{duality} that the Shi bijection induces the bijection  
between $\AD_0$ and the bounded dominant regions. 
We survey these results in Section~\ref{old}.
\\[.5ex]
The goal of this paper is two-fold. First, we study the
ideals lying inside $\gH$, the set
of all positive roots that are not orthogonal to $\theta$.
Such ideals are said to be {\it of Heisenberg type\/}.
Second, we study the {\it minimax\/} elements of $\HW$, i.e., those
that are simultaneously minimal and maximal.
The corresponding strictly positive ideals
are said to be minimax, too. 
Both minimal and maximal elements of $\HW$
are particular instances of {\it dominant\/} elements. To any dominant 
element one may attach an ideal, and 
we consider dominant elements of special sort,
which are said to be of {\it Heisenberg type\/},
see Definition~\ref{def-heis}. The corresponding ideals
lie in $\gH$, hence the term. 
We obtain a explicit description of dominant elements
of Heisenberg type and distinguish minimal and maximal ones among them.
Our description says, in particular, that the dominant (resp.
minimal) elements of Heisenberg type 
are parameterized by $\Delta_l$ (resp. $\Delta_l\setminus (-\Pi)$), 
where $\Delta_l$ is the set of long roots.
We also prove that the minimal elements of Heisenberg type are in
a bijection with the non-trivial ideals of Heisenberg type and give 
explicit formulae for these ideals, see Theorem~\ref{main3}.
It follows that the number of non-trivial ideals
inside $\gH$ is equal to $\#(\Delta_l\setminus \Pi)$.
\\[.5ex]
In Section~\ref{minm}, combining known results for minimal
and maximal elements, we obtain a characterization of minimax elements
(ideals),
and establish a bijection between the minimax elements of $\HW$
and the the coroot lattice points lying in a polytope $D_{mm}$.
If $V$ is the real vector space spanned by $\Delta$ in $\te^*$, then
\[
  D_{mm}:=\{x\in V \mid  -1\le (\ap,x)\le 1\ \ \forall \ap\in\Pi
  \quad\&\quad  0\le (\theta,x)\le 2 \} 
\]
It is the intersection of two simplices corresponding to the
minimal and maximal elements of $\HW$.
The geometric meaning of minimax ideals is that, under the Shi bijection, 
they correspond to the regions consisting of a single alcove.
It would be interesting to further investigate the configuration
of those dominant regions that consist of a single alcove.
We also give an upper bound on the sum $\#\Gamma(I)+\#\Xi(I)$, where
$\Gamma(I)$ is the set of generators (minimal elements) of $I$ and
$\Xi(I)$ is the set of maximal elements of $\Delta^+\setminus I$.
Since $\Pi$ is the only antichain of cardinality $\rk\g$, one 
readily sees that
the maximal possible value of the above sum is $2\rk\g-1$, which is attained
for $I=\Delta^+\setminus\Pi$. However, for the minimax
ideals this sum is at most $\rk\g+1$
(at most $\rk\g$, if $\g\ne {\frak sl}_{2n+1}$), see Proposition~\ref{bound-mm}.
\\[.5ex]
Write  $\AD_{mm}$ or $\AD_{mm}(\g)$ for the set of minimax ideals of $\g$.
In Section~\ref{count}, we compute  $\#\AD_{mm}(\g)$ for
all simple Lie algebras. 
Using the coefficients of $\theta$, one may form a Laurent polynomial (in $x$).
We prove that $\#\AD_{mm}$ equals the coefficient of $x$ in this polynomial,
divided by the index of connection of $\Delta$.
It turns out that some famous integer sequences emerge in connection with 
the minimax ideals in classical Lie algebras. Namely,

$\#\AD_{mm}({\frak sl}_{n+1})=M_n$, the
$n$-{\it th Motzkin number\/};  

$\#\AD_{mm}({\frak sp}_{2n})= 
\#\AD({\frak so}_{2n+1})=\text{dir}_n$, the number
of {\it directed animals of size\/} $n$.
\\
We also show  that $\#\AD_{mm}({\frak so}_{2n})=2\text{dir}_{n-2}+
\text{dir}_{n-1}$.
See \cite{martin},\,\cite{dsh},\,\cite{anim}, and \cite[Ch.\,6]{rstan}
for relevant background and numerous combinatorial interpretations of 
these numbers. 
In Section~\ref{sln}, an explicit matrix characterization of 
the minimax ideals for ${\frak sl}_{n+1}$ and ${\frak sp}_{2n}$ is
obtained. The description for ${\frak sl}_{n+1}$ can be stated as follows.
Let $\{\ap_i=\esi_i-\esi_{i+1}\mid i=1,\dots,n\}$ 
be the standard set of simple roots.
Let $\Gamma=\{\gamma_1,\dots,\gamma_k\}$ be the set of generators of an 
ideal $I$, where $\gamma_t=\ap_{i_t}+\ap_{i_t+1}+
\ldots +\ap_{j_t}$, $t=1,2,\dots,k$.
Then $I$ is minimax if and only if $j_t\ne i_{s}$ for all pairs $(t,s)$.
In case $t=s$, this means that  $\gamma_t$ cannot be a simple root. 
Using these descriptions, we also compute the generating function 
for the statistic ``the number of generators'' on $\AD_{mm}(\g)$ for
$\g={\frak sl}_{n+1}$ and ${\frak sp}_{2n}$.

{\small {\bf Acknowledgements.} A part of this paper was written during 
my stay at the Max-Planck-Institut f\"ur Mathematik (Bonn). 
I would like to thank this institution for hospitality
and excellent working conditions. 
}


\section{Notation and other preliminaries}
\label{prelim}

\noindent
\begin{subs}{Main notation}
\end{subs}
$\Delta$ is the root system of $(\g,\te)$ and
$W$ is the usual Weyl group. For $\ap\in\Delta$, $\g_\ap$ is the
corresponding root space in $\g$.

$\Delta^+$  is the set of positive
roots and $\rho=\frac{1}{2}\sum_{\ap\in\Delta^+}\ap$.

$\Pi=\{\ap_1,\dots,\ap_p\}$ is the set of simple roots in $\Delta^+$
and $\theta$ is the highest root in $\Delta^+$.

$e_1,e_2,\dots,e_p$ \ are the exponents and $h$ is the Coxeter
number of $W$.
 \\
 We set $V:=\te_{\Bbb R}=\oplus_{i=1}^p{\Bbb R}\ap_i$ and denote by
$(\ ,\ )$ a $W$-invariant inner product on $V$. As usual,
$\mu^\vee=2\mu/(\mu,\mu)$ is the coroot
for $\mu\in \Delta$.

${\gC}=\{x\in V\mid (x,\ap)>0 \ \ \forall \ap\in\Pi\}$
\ is the (open) fundamental Weyl chamber.

${\gA}=\{x\in V\mid (x,\ap)>0 \ \  \forall \ap\in\Pi \ \ \& \ 
(x,\theta)<1\}$ \ is the fundamental alcove.


$Q^+=\{\sum_{i=1}^p n_i\ap_i \mid n_i=0,1,2,\dots \}$
and $Q^\vee=\oplus _{i=1}^p {\Bbb Z}\ap_i^\vee\subset V$  
is the coroot lattice.
\\
Letting $\widehat V=V\oplus {\Bbb R}\delta\oplus {\Bbb R}\lb$, we extend
the inner product $(\ ,\ )$ on $\widehat V$ so that $(\delta,V)=(\lb,V)=
(\delta,\delta)=
(\lb,\lb)=0$ and $(\delta,\lb)=1$.

$\widehat\Delta=\{\Delta+k\delta \mid k\in {\Bbb Z}\}$ is the set of affine
real roots and $\widehat W$ is the  affine Weyl group.
\\
Then $\widehat\Delta^+= \Delta^+ \cup \{ \Delta +k\delta \mid k\ge 1\}$ is
the set of positive
affine roots and $\widehat \Pi=\Pi\cup\{\ap_0\}$ is the corresponding set
of affine simple roots,
where $\ap_0=\delta-\theta$.
The inner product $(\ ,\ )$ on $\widehat V$ is
$\widehat W$-invariant. The notation $\beta>0$ (resp. $\beta <0$)
is a shorthand for $\beta\in\HD^+$ (resp. $\beta\in -\HD^+$).
\\
For $\ap_i$ ($0\le i\le p$), we let $s_i$ denote the corresponding simple
reflection in $\widehat W$.
If the index of $\ap\in\widehat\Pi$ is not specified, then we merely write
$s_\ap$. 
The length function on $\widehat W$ with respect
to  $s_0,s_1,\dots,s_p$ is denoted by $\ell$.
For any $w\in\widehat W$, we set
\[
   N(w)=\{\ap\in\widehat\Delta^+ \mid w(\ap) \in -\widehat \Delta^+ \} .
\]
It is standard that $\#N(w)=\ell(w)$ and $N(w)$ is {\it bi-convex\/}. 
The latter means
that both $N(w)$ and $\HD^+\setminus N(w)$ are subsets of $\HD^+$
that are closed under addition.
Furthermore, the assignment $w\mapsto N(w)$ sets up
a bijection between the elements of $\HW$ and the finite bi-convex subsets
of $\HD^+$.

\begin{subs}{Ideals and antichains}               \label{ideals}
\end{subs}
Recall that $\be$ is the Borel subalgebra of $\g$ corresponding 
to $\Delta^+$
and $\ut=[\be,\be]$. 
Let $\ce\subset\be$ be an {\sf ad}-nilpotent ideal. Then
$\ce=\underset{\ap\in I}{\oplus}\g_\ap$
for some $I\subset \Delta^+$. This $I$ is said to be an {\it ideal\/} (of
$\Delta^+$). More precisely, a set $I\subset\Delta^+$ is an ideal, if
whenever $\gamma\in I,\mu\in\Delta^+$, and $\gamma+\mu\in\Delta$, then
$\gamma+\mu\in I$.
Our exposition will be mostly combinatorial, i.e., in place of
{\sf ad}-nilpotent ideal of $\be$ we will deal with the respective 
ideals of $\Delta^+$.
\\[.5ex]
For $\mu,\gamma\in\Delta^+$, write $\mu\curle\gamma$, if
$\gamma-\mu\in Q^+$. The notation $\mu\prec\gamma$ means that
$\mu\curle\gamma$ and $\gamma\ne\mu$.
We regard $\Delta^+$ as poset under `$\cyeq$'.
Let $I\subset\Delta^+$ be an ideal. An element $\gamma\in I$ is called
a {\it generator\/}, if $\gamma-\ap\not\in I$ for any $\ap\in\Pi$.
In other words, $\gamma$ is a minimal element of $I$.
We write $\Gamma(I)$ for the set of generators of $I$.
It is easily seen that $\Gamma(I)$ is an {\it antichain\/} 
of $\Delta^+$, i.e., $\gamma_i\not\curle\gamma_j$ for any pair 
$(\gamma_i,\gamma_j)$ in $\Gamma(I)$.
Conversely, if $\Gamma\subset \Delta^+$ is an antichain, 
then the ideal 
\[
I\langle\Gamma\rangle:=\{\mu\in \Delta^+\mid \mu\succcurlyeq\gamma_i
\text{ for some } \gamma_i\in \Gamma \} 
\]
has $\Gamma$ as the set of generators.
Let $\AN$ denote the set of all antichains in $\Delta^+$.
In view of the above bijection  
$\AD \overset{1:1}{\longleftrightarrow}\AN$, 
we will freely switch between
ideals and antichains. An ideal $I$ is called {\it strictly positive\/}, if
$I\cap\Pi=\varnothing$. The set of strictly positive ideals is denoted
by $\AD_0$.

Given $I\in\AD$, write $\Xi(I)$ for the set of maximal elements of
$\Delta^+\setminus I$.


\section{Minimal and maximal elements in affine Weyl groups and
dominant regions of the Catalan arrangement}  
\label{old}
\setcounter{equation}{0}

\noindent
This section is a review of known results.
As is well known, $\widehat W$ is isomorphic to a semi-direct product
of $W$ and $Q^\vee$ \cite{hump}.
Given $w\in\widehat W$, there is a unique decomposition
\begin{equation}  \label{affine}
w=v{\cdot}t_{r}\ ,
\end{equation}
where $v\in W$ and  $t_{r}$ is the translation
corresponding to $r\in Q^\vee$. Then $w^{-1}=v^{-1}{\cdot}t_{-v(r)}$.
The word ``translation" means the following.  The group $\HW$ has two 
natural actions:

(a) the linear action on $\widehat V=
V\oplus{\Bbb R}\delta\oplus{\Bbb R}\lb$;

(b) the affine-linear action on $V$.
\\
We use `$\ast$' to denote the second action.
For the linear action, we have 
$w^{-1}(x)=v^{-1}(x)+(x, v(r))\delta$ \ for any
$x\in V\oplus {\Bbb R}\delta$. In particular,
\begin{equation} \label{simroots} \begin{array}{l}
  w^{-1}(\ap_i)=v^{-1}(\ap_i)+(\ap_i, v(r))\delta, \quad i\ge 1 ,\\
  w^{-1}(\ap_0)=-v^{-1}(\theta)+(1-(\theta, v(r)))\delta \ .
\end{array}
\end{equation}
While the affine-linear action is 
given by $w^{-1}\asst x=v^{-1}(x)-r$ for $x\in V$.  In particular,
$t_r\asst y=y+r$, so that $t_r$ is a true translation for the
$\ast$-action on $V$. 
\\[.5ex]
Let us say that $w\in\HW$ is {\it dominant\/}, if
$w(\ap)>0$ for all $\ap\in\Pi$. 
Obviously, $w$ is dominant if and only if
$N(w)\subset \cup_{k\ge 1} (k\delta-\Delta^+)$.
It also follows from \cite[1.1]{CP1} that
$w$ is dominant if and only if 
$w^{-1}\asst \gA\subset \gC$. Write $\HW_{dom}$ for the set of
dominant elements.

\begin{s}{Proposition {\ququ \cite{losh}}}  \label{opis-dom}
\begin{itemize}
\item[\sf (i)] If $w=v{\cdot}t_r\in \HW_{dom}$, then $r\in -\ov{\gC}$;
\item[\sf (ii)] The mapping $\HW_{dom}\to Q^\vee$ given
by $w=v{\cdot}t_r\mapsto v(r)$ is a bijection.
\end{itemize}
\end{s}%
Letting $\delta-I:=N(w)\cap (\delta-\Delta^+)$, we easily
deduce that $I$ is an ideal, if $w\in \HW_{dom}$. We say $\delta-I$ is 
the {\it first layer\/} of $N(w)$ and $I=I_w$ is the 
{\it first layer ideal\/} of $w$.
However, an ideal $I$ may well arise from different dominant elements.
To obtain a bijection, one has to impose further constraints on
dominant elements. One may attempt to consider either maximal or minimal
bi-convex subsets with first layer $\delta-I$.
This naturally leads to notions of
`minimal' and `maximal' elements. This terminology suggested in
\cite{eric} is also explained by a relationship between these elements 
and dominant regions of the Shi arrangement, see Theorem~\ref{min+max}.

\begin{rem}{Definition}   \label{def-mi}
$w\in\HW$ is called {\it minimal\/}, if
\begin{itemize}
\item[\sf (i)] \ $w$ is dominant;
\item[\sf (ii)] \ if $\ap\in\HP$ and  
$w^{-1}(\ap)=k\delta+\mu$ for some $\mu\in \Delta$, then $k\ge -1$.
\end{itemize}
\end{rem}%
Using (i), condition (ii) can be made more precise. If $k\in\{-1,0\}$,
then $\mu\in \Delta^+$.
\\
The set of minimal elements is denoted by $\HW_{min}$.
For any $\gamma\in I$, define the number
$l(\gamma,I)$ as follows
\[
 l(\gamma,I)=\max\{ m\mid \gamma=\sum_{i=1}^m \varkappa_i,\ \ 
   \varkappa_i\in I\} \ .
\]
\vskip-1ex
\begin{s}{Proposition {\ququ \cite[Prop.\,2.12]{CP1}}}
\label{bij-mi}
There is a bijection between $\HW_{min}$ and\/ $\AD$.
Namely,
\begin{itemize}
\item given $w\in \HW_{min}$, the corresponding ideal is
$\{\mu\in\Delta^+\mid \delta-\mu\in N(w)\}$;
\item given $I\in\AD$, the corresponding minimal element is determined by
the finite bi-convex set 
\[
   \bigcup_{k\ge 1}(k\delta-I^k)=
\{ m\delta-\gamma\mid \gamma\in I \ \ \& \ \ 
 1\le m\le l(\gamma,I) \} \subset \HD^+ \ .
\] 
Here $I^k$ is defined inductively by $I^k=(I^{k-1}+I)\cap\Delta^+$.
\end{itemize}
\end{s}%
The minimal element corresponding to $I\in\AD$ is denoted by $w_{min}(I)$.
Conversely, the first layer ideal of $w\in\HW_{min}$ (= the ideal
corresponding to $w$) is denoted by $I_w$.
\\
If $N\subset \HD^+$ is an arbitrary finite convex subset, containing 
$\delta-I$, then it must also contain 
$N(w_{min}(I))=\cup_{i\ge 1}(k\delta-I^k)$. So, the 
latter is the minimal bi-convex subset containing $\delta-I$.
\\[.5ex]
We have also the following
\begin{s}{Corollary}  \label{ell}
For $w\in \HW_{min}$, we have $\ell(w)=\sum_{k\ge 1} \#((I_w)^k)$.
\end{s}%
In terms of minimal elements, one can give an explicit description 
of the generators of an ideal.

\begin{s}{Proposition {\ququ \cite[Theorem\,2.2]{duality} 
\cite[Cor.\,6.3(1)]{eric}}} 
\label{gen-mi}  \par
If $w\in\HW_{min}$, then
$\Gamma(I_w)=\{\gamma\in\Delta^+ \mid w(\delta-\gamma) \in-\HP\}$.
\end{s}%
A geometric description of the minimal elements relates
them to the integral points in a certain simplex.

\begin{s}{Proposition {\ququ \cite[Prop.\,2 \& 3]{CP2}}} 
\label{opis-mi} 
Set $D_{min}=\{x\in V \mid (x,\ap)\ge -1 \ \forall\ap\in\Pi \ \ \& \ 
\ (x,\theta)\le 2\}$. Then \par
1. $w=v{\cdot}t_r\in \HW_{min}$ \ $\Longleftrightarrow
\left\{ \begin{array}{l} w \text{ is dominant},  \\
                         v(r)\in D_{min}\cap {Q}^\vee \ .
       \end{array}\right.$

2. The mapping $\HW_{min} \to D_{min}\cap Q^\vee$, 
$w=v{\cdot}t_r\mapsto v(r)$, is a bijection.
\end{s}%
It follows that $\#(\AD)$ equals the number
of integral points in $D_{min}$. (Unless otherwise stated, an
'integral point' is a point lying in $Q^\vee$.)
A pleasant feature of this situation is that there
is an element of $\HW$ that takes $D_{min}$ to a dilated
closed fundamental alcove. Namely,
$w(D_{min})=(h+1)\ov{\gA}$ for some $w\in \HW$, see \cite[Thm.\,1]{CP2}.
Write $\theta$ as a linear combination of simple roots:
$\theta= \sum_i c_i\ap_i$. The integers $c_i$ are said to be the
{\it coordinates of\/} $\theta$. 
By a result of M.\,Haiman \cite[7.4]{mark}, the number of integral points
in $t\ov{\gA}$ is equal to 
\begin{equation}  \label{haiman}
\prod_{i=1}^p \frac{t+e_i}{1+e_i} \ 
\end{equation}
whenever $t$ is relatively prime with all the coordinates 
of $\theta$. Since this condition is satisfied for
$t=h+1$, one obtains, see \cite{CP2},
\begin{equation}  \label{chislo-mi}
   \#\HW_{min}=\#\AD=\prod_{i=1}^p \frac{h+e_i+1}{e_i+1} \ .
\end{equation}
Maximal elements of $\HW$ are introduced by E.\,Sommers in \cite{eric}.
He also obtained  main results for these elements. 

\begin{rem}{Definition}   \label{def-ma}
$w\in\HW$ is called {\it maximal\/}, if
\begin{itemize}
\item[\sf (i)] $w$ is dominant;
\item[\sf (ii)] if $\ap\in\HP$ and  
$w^{-1}(\ap)=k\delta+\mu$ for some $\mu\in \Delta$, then $k\le 1$.
\end{itemize}
\end{rem}%
Using (i), condition (ii) can be made more precise. 
If $k=1$, then $\mu\in -\Delta^+$; if $k=0$, then $\mu\in \Delta^+$
The set of maximal elements is denoted by
$\HW_{max}$.

If $I\in\AD_0$, then for any $\mu\in\Delta^+$ we define the number
$k(\mu, I)$ as follows:
\[
 k(\mu,I)=\min\{ n\mid \mu=\sum_{i=1}^n \nu_i,\quad
   \nu_i\in \Delta^+\setminus I\} \ .
\]
Notice that this definition makes sense only for strictly positive ideals.

\begin{s}{Proposition {\ququ \cite[Sect.\,5]{eric}}}
\label{bij-ma}
There is a bijection between $\HW_{max}$ and $\AD_0$.
Namely,
\begin{itemize}
\item \ given $w\in \HW_{min}$, the corresponding strictly
positive ideal is
$\{\mu\in\Delta^+\mid \delta-\mu\in N(w)\}$;
\item \ given $I\in\AD_0$, the corresponding maximal element is 
determined by the finite bi-convex set 
\end{itemize}
\hbox to \textwidth{\enspace ($\Diamond$)\hfil
 $\{ m\delta-\gamma\mid \gamma\in I \quad\&\quad 
 1\le m\le k(\gamma,I){-}1\}$. \hfil}
\end{s}%
The maximal element corresponding to $I\in\AD_0$ is denoted by 
$w_{max}(I)$. Accordingly, the strictly positive ideal corresponding to 
$w\in\HW_{max}$ (the first layer ideal of $w$) is  $I_w$. 
It follows from the previous discussion that to any $I\in\AD_0$ one 
can assign two elements of $\HW$; namely, 
$w_{min}(I)$ and $w_{max}(I)$.
\\[.5ex]
Let $w\in\HW$ be any dominant element with
first layer ideal $I\in\AD_0$. Since $\HD^+\setminus N(w)$ is convex and
contains $\delta-(\Delta^+\setminus I)$, it follows 
from the very definition of numbers 
$k(\gamma,I)$ that $m\delta-\gamma\in \HD^+\setminus N(w)$
for all $m\ge k(\gamma,I)$. Hence $N(w)$ is contained in the finite
set given by Eq.~\ref{bij-ma}($\Diamond$).
This shows that $N(w_{max}(I))$ is the maximal bi-convex subset
with first layer $\delta-I$. Remember also that 
$N(w_{min}(I))$ is the minimal bi-convex subset
with first layer $\delta-I$. Hence
$N(w_{min}(I))\subset N(w_{max}(I))$ and therefore
\begin{equation}  \label{l<k}
    l(\gamma,I)\le k(\gamma,I)-1 \quad \text{for any\ }\ \gamma\in I\ .
\end{equation}
Recall that $\Xi(I)$ is the set of maximal elements of $\Delta^+\setminus
I$. If $I\in\AD_0$, then a description of $\Xi(I)$ can be given in terms
of $w_{max}(I)$:

\begin{s}{Proposition {\ququ \cite[6.3(2)]{eric}}}  \label{gen-ma}
If $w\in\HW_{max}$, then
$\Xi(I_w)=\{\gamma\in\Delta^+\mid w(\delta-\gamma)\in\HP\}$.
\end{s}%
Next, we recall a geometric characterization of the maximal elements.

\begin{s}{Proposition {\ququ (cf. \cite[Prop.\,5.6]{eric})}} 
\label{opis-ma} 
Set $D_{max}=\{x\in V \mid (x,\ap)\le 1 \ \forall\ap\in\Pi \ \ \& \ \ 
(x,\theta)\ge 0\}$. Then \par
1. $w=v{\cdot}t_r\in \HW_{max}$ \ $\Longleftrightarrow
\left\{ \begin{array}{l} w \text{ is dominant},  \\
                         v(r)\in D_{max}\cap {Q}^\vee \ .
       \end{array}\right.$

2. The mapping $\HW_{max} \to D_{max}\cap Q^\vee$, 
$w=v{\cdot}t_r\mapsto v(r)$, is a bijection.
\end{s}%
In order to compute  $\#(D_{max}\cap Q^\vee)$,  
we replace $D_{max}$ with another simplex.
Let $\{\varpi^\vee_i\}_{i=1}^p$ denote the dual basis of $V$ for 
$\{\ap_i\}_{i=1}^p$. Set $\rho^\vee=\sum_{i=1}^p\varpi^\vee_i$.
Since the sum of the coordinates of $\theta$ equals $h-1$, the translation 
$x\mapsto t_{-\rho^\vee}\asst x=x-\sum_{i=1}^p\varpi^\vee_i$ 
takes $D_{max}$ to the negative dilated fundamental alcove
\[
-(h-1)\ov{\gA}=\{ x\in V\mid (x,\ap)\le 0 \quad \forall\,\ap\in\Pi;\ 
   (x,\theta)\ge 1-h\} .
\]
It may happen that $\rho^\vee$ does not belong to $Q^\vee$, so that
this translation, which is in
the extended affine Weyl group, does not belong to $\HW$,
while we wish to have a transformation from $\HW$.
Nevertheless, since $h-1$ is relatively prime with the index of
connection of $\Delta$, it follows from \cite[Lemma\,1]{CP2} that 
there is an element of
$\HW$ that takes $D_{max}$ to $(1-h)\ov{\gA}$.

Using again the above-mentioned result of Haiman, see Eq.~\re{haiman},
one obtains, see \cite{ath02},\,\cite{duality},\,\cite{eric},   

\begin{equation} \label{chislo-ma}
   \#\HW_{max}=\#(\AD_0)= \prod_{i=1}^p \frac{h+e_i-1}{e_i+1}  \ .
\end{equation}
Recall a bijection between $\AD$
and the dominant regions of the Catalan arrangement,
which is due to 
Shi \cite[Theorem\,1.4]{shi}.  
\\[.5ex]
For $\mu\in\Delta^+$ and $k\in{\Bbb Z}$, define the hyperplane
$\cH_{\mu,k}$ in $V$ as $\{x\in V \mid (x,\mu)=k\}$.
The {\it Catalan arrangement\/}, $\mathrm{Cat}(\Delta)$,
is the collection of hyperplanes 
$\cH_{\mu,k}$, where $\mu\in \Delta^+$ and $k=-1,0,1$.
The {\it regions\/} of an arrangement are the connected components of
the complement in $V$ of the union of all its hyperplanes. 
Any region lying in $\gC$ is said to be
{\it dominant\/}.
Obviously, the dominant regions of $\mathrm{Cat}(\Delta)$ are the same as 
those for the {\it Shi arrangement\/} $\mathrm{Shi}(\Delta)$.
The latter is the collection of hyperplanes
$\cH_{\mu,k}$, where $\mu\in \Delta^+$ and $k=0,1$.
But, it is sometimes more convenient to deal with the arrangement 
$\mathrm{Cat}(\Delta)$, since it is $W$-invariant.
\\[.5ex]
The Shi bijection takes an ideal $I\subset \Delta^+$ to the dominant region
\begin{equation}   \label{bij-shi}
 R_I =\{ x\in \gC \mid  (x,\gamma)>1, \text{ if \ }\gamma\in I \quad \&\quad
           (x,\gamma)<1,  \text{ if \ }\gamma\not\in I
           \}\ .
\end{equation}
A region of an arrangement is called {\it bounded\/}, if it is contained 
in a sphere about the origin.

\begin{s}{Proposition {\ququ \cite{ath02},\cite{duality}}}  
\label{regions}
$I\in\AD(\g)_0$ 
if and only if the region $R_I$ is bounded.
\end{s}%
A relationship between the theory of minimal/maximal elements and
regions of the Catalan arrangement is as follows.

\begin{s}{Theorem} \label{min+max} \\
(i) \ \cite{CP2} \ Suppose $I\in\AD$ and $w=w_{min}(I)$. 
Then $w^{-1}\asst\gA$ is the alcove closest
to the origin in $R_I$;
\\[.6ex]
(ii) \ \cite{eric} 
Suppose $I\in\AD_0$ and $w=w_{max}(I)$. Then $w^{-1}\asst\gA$ is the alcove
most distant from the origin in $R_I$.
\end{s}%
%


\section{Dominant elements and ideals of Heisenberg type}
\label{heis}
\setcounter{equation}{0}

\noindent
Given $w\in\HW$, write $w(\ap_0)=-m\delta+\nu$, where $m\in {\Bbb Z}$
and $\nu\in \Delta$. Since $(\theta,\theta)=(\ap_0,\ap_0)=(\nu,\nu)$, 
the root $\nu$ is long.
The root $\nu$ is called the {\it rootlet\/} of $w$, denoted ${\rt}(w)$.
In what follows, $\Delta_l$ stands for the set of long roots.
In general, the rootlet can be any root in $\Delta_l$; 
but, for various classes of elements of $\HW$, different constraints emerge.

\begin{s}{Lemma}   \label{koreshki}
\begin{itemize}
\item[\sf (i)] \ If $1\ne w\in\HW_{dom}$, then $m\ge 1$;
\item[\sf (ii)] \ If $w\in\HW_{min}$, then ${\rt}(w)\not\in -\Pi$;
\item[\sf (iii)] \  If $w\in\HW_{max}$, then ${\rt}(w)\not\in \Pi$.
\end{itemize}
\end{s}\begin{proof}
(i) \ Recall the ``affine-linear decomposition"
$w=\bar w{\cdot}t_r$, where $\bar w\in W$ and 
$r\in -\ov{\gC}$. If $w\ne 1$, then $r\ne 0$, see Prop.~\ref{opis-dom}.
Hence $(\theta, r)<0$. Furthermore, $-\ov{\gC}$ contains no elements
$x\in Q^\vee$ such that $(x,\theta)=-1$. For, Haiman's 
formula~\re{haiman} shows
that $\ov{\gA}$ contains a unique element of $Q^\vee$, namely, the origin.
Hence $(r,\theta)\le -2$. It remains to observe that
$w(\ap_0)=-\bar w(\theta)+(1+(\theta,r))\delta$. Thus,
$-m=1+(\theta,r)\le -1$.
\\[.5ex]
(ii) If $w(\ap_0)=-m\delta-\ap$, $\ap\in\Pi$, then
$w^{-1}(\ap)=-(m+1)\delta+\theta$ and $-(m+1)\le -2$.
Hence $w\not\in \HW_{min}$.
\\[.5ex]
(iii) If $w(\ap_0)=-m\delta+\ap$, $\ap\in\Pi$, then
$w^{-1}(\ap)=(m+1)\delta-\theta$ and $m+1\ge 2$.
Hence $w\not\in \HW_{max}$.
\end{proof}%
{\bf Remark.} Notice that the formulas for $w(\ap_0)$ shows that
${\rt}(w)=-\bar w(\theta)$.
\\[.5ex]
Because there is a bijection between the minimal 
elements of $\HW$ and the ideals of $\Delta^+$, 
one may define the rootlet for any ideal. That is,
given $w\in\HW_{min}$, we set ${\rt}(I_w):={\rt}(w)$. It should be
noticed that if $I\in\AD_0$, then the roots ${\rt}(w_{min}(I))$ and
${\rt}(w_{max}(I))$ are usually different. So, if we wish to get an
unambiguous notion of the rootlet for any ideal, then the minimal elements
have to be exploited.
\\[.5ex]
Amongst all ideals of $\Delta^+$, we distinguish the ideal consisting of
all roots that are not orthogonal to $\theta$. Set
\[
\gH=\{\gamma\in\Delta^+\mid (\gamma,\theta)>0\} \ .
\]
The corresponding subspace of $[\be,\be]$ is the standard Heisenberg 
subalgebra, so that we say that $\gH$ is the {\it Heisenberg ideal\/}.
Notice that $\gH^2=\{\theta\}$ and $\gH^3=\varnothing$.

\begin{rem}{Definition}   \label{def-heis}
We say that $w\in\HW$ is {\it of Heisenberg type\/}, if $w=vs_0$ for
some $v\in W$. An ideal is said to be {\it of Heisenberg type\/}, if
it is contained in $\gH$.
\end{rem}%
The term for elements of $\HW$ is explained by the following

\begin{s}{Lemma}
If $w\in \HW_{dom}$ is of Heisenberg type, then the first layer ideal
of $w$ lies in $\gH$.
\end{s}\begin{proof}
Recall that the first layer ideal of $w$ consists of all 
$\gamma\in \Delta^+$ such that $w(\delta-\gamma)<0$. 
If $w$ is of the form $vs_0$ \ ($v\in W$) and
$(\gamma,\theta)=0$, then $w(\delta-\gamma)=\delta- v(\gamma)>0$.
\end{proof}%
The main goal of this section is to give a characterization of the
dominant elements of Heisenberg type, and then to describe all ideals 
in $\gH$. 

\begin{s}{Proposition}  \label{heis0}
Let $w=vs_0\in \HW$ be of Heisenberg type. Then
\begin{itemize}
\item[\sf (i)] \ $w\in \HW_{dom}$ \ $\Leftrightarrow$ \ $v(\ap)>0$ for any
$\ap\in\Pi$ such that $(\ap,\theta)=0$;
\item[\sf (ii)] \ ${\rt}(w)=v(\theta)$;
\item[\sf (iii)] \ $w\in \HW_{min}$ \ $\Leftrightarrow$ \ $w\in \HW_{dom}$
and $v(\theta)\not\in -\Pi$;
\item[\sf (iv)] \ $w\in \HW_{max}$ \ $\Leftrightarrow$ \ $w\in \HW_{dom}$,
$(v(\theta),\theta)\ge 0$, and $v(\theta)\not\in \Pi$.
\end{itemize}
\end{s}\begin{proof}
(i) The property of being dominant means $w(\ap)>0$ for all $\ap\in\Pi$.
If $(\ap,\theta)>0$, then $w(\ap)=v(\delta-(\theta-\ap))=\delta-
v(\theta-\ap)>0$, i.e., it is always satisfied.
If $(\ap,\theta)=0$, then $w(\ap)=v(\ap)$. Hence the condition.

(ii) \ $vs_0(\ap_0)=v(\theta)-\delta$.

(iii) \& (iv). The conditions of minimality and maximality impose 
constraints on the coefficient of $\delta$ in $w^{-1}(\ap)$, $\ap\in\HP$.
\\[.5ex]
For $\ap\in\Pi$, we have $w^{-1}(\ap)=s_0{\cdot}v^{-1}(\ap)$.
If $v^{-1}(\ap)=\pm\theta$, then $s_0(\pm\theta)=\pm(2\delta-\theta)$, which
is bad. More precisely, $w$ is not maximal, if $v(\theta)=\ap$;
and $w$ is not minimal, if  $v(\theta)=-\ap$. These are the only bad 
possibilities.
\\[.5ex]
For $\ap_0$, we have $w^{-1}(\ap_0)=s_0(\delta-v^{-1}(\theta))
=\delta-v^{-1}(\theta)-(\theta,v(\theta)^\vee)(\delta-\theta)$.
It is easily seen that the only bad possibility is $(\theta,v(\theta))<0$,
in which case $w$ fails to be maximal.
\end{proof}%
Thus, the conditions of maximality and minimality for the Heisenberg type
elements are stated in terms of the rootlet.
\\[.6ex]
Recall the following standard fact: \\[.6ex]
\refstepcounter{equation} (\theequation) \label{uv}
\qquad
\parbox{410pt}{{\it Suppose $u,v\in W$. Then $\ell(u)+\ell(v)=\ell(uv)$
if and only if $N(u)\cap N(v^{-1})=\varnothing$. If these conditions are
satisfied, then $N(uv)=N(v)\sqcup v^{-1}(N(u))$.}} 
\\[.6ex]
If $\ell(u)+\ell(v)=\ell(uv)$, then we say that $uv$ is a {\it reduced 
decomposition\/} for this product.
We will also repeatedly use the following result, see
\cite[Theorem\,4.1]{lp}:
\\[.6ex]
\refstepcounter{equation} (\theequation) \label{-1}
\qquad
\parbox{420pt}
{{\it Given $\nu\in \Delta^+_l$, there is a unique shortest element 
in $W$ taking $\theta$ to $\nu$, denoted by $w_\nu$.
We have $N(w_\nu^{-1})=\{\gamma\in\Delta^+\mid (\gamma,\nu^\vee)=-1 \}$. }}
\\[.6ex]
Given $\nu\in\Delta^+$, write $s_\nu$ for the corresponding reflection
in $W$. 

\begin{s}{Lemma}    \label{dlina}
Suppose $\nu\in\Delta^+_l$. Then
\begin{itemize} 
\item[\sf (i)] \ $\ell(s_\nu)=2(\rho,\nu^\vee)-1$;
\item[\sf (ii)] \ $\ell(s_\nu w_\nu)=\ell(s_\nu)+\ell(w_\nu)=
(\rho,\theta^\vee+\nu^\vee)-1$.
\item[\sf (iii)] \ $N(s_\nu)=\{\nu\}\cup
\{\gamma\in\Delta^+\mid (\gamma,\nu^\vee)=1 \ \ 
\& \ \ \gamma \prec \nu\}$.
\end{itemize}
\end{s}\begin{proof}
(i) We have $s_\nu(-\nu)=\nu$.
Recall that $(\rho,\ap^\vee)=1$ for all $\ap\in\Pi$.
We will argue by induction on $(\rho,\nu^\vee)$. If $\nu\in\Pi$, 
then $\ell(s_\nu)=1$, and the claim is true. 
\\
A straightforward calculation shows that
$(\rho, s_\ap(\mu)^\vee)= (\rho,\mu^\vee)-(\ap, \mu^\vee)$
for $\mu\in\Delta$ and $\ap\in\Pi$.
If $\mu$ is long, we have
$(\ap, \mu^\vee)\ge -1$ unless $\mu=-\ap$.
Hence $(\rho, s_\ap(\mu)^\vee)\le (\rho,\mu^\vee)+1$, if 
$\mu\ne -\ap$, and $(\rho, s_\ap(-\ap)^\vee)=(\rho,-\ap^\vee)+2$.
Let $s_\nu=s_{i_1}\ldots s_{i_k}$ be a reduced decomposition.
Consider the corresponding sequence of roots:
\[
  (-\nu) \stackrel{s_{i_k}}{\to}{\circ} \stackrel{s_{i_{k-1}}}{\to}{\circ}
\ldots \stackrel{s_{i_{2}}}{\to}\circ\stackrel{s_{i_1}}{\to} \nu \ .
\]
The level function $\mu\mapsto (\rho,\mu^\vee)$ attains $2(\rho,\nu^\vee)$
values between $-(\rho,\nu^\vee)$ and $(\rho,\nu^\vee)$, 
because the zero level is missing.
Moreover, it follows from the previous inequalities that at each step
one can go up at most to the next existing level.
Hence, $\ell(s_\nu)\ge 2(\rho,\nu^\vee)-1$. 
On the other hand, $s_\nu=s_\ap s_{\gamma}s_\ap$, where 
$\gamma=s_\ap(\nu)$. Because $\nu\not\in\Pi$, one can 
find an $\ap\in\Pi$ so that $(\ap,\nu^\vee)=1$. Then
we obtain $(\rho,\gamma^\vee)=(\rho,\nu^\vee)-1$ and
by the induction assumption 
$\ell(s_\nu)\le 2+\ell(s_\gamma)=2(\rho,\nu^\vee)-1$.
This completes the inductive step.

(ii) The first equality follows from Eq.~\re{uv}. Indeed,
$N(w^{-1}_\nu)=\{\gamma\mid (\gamma,\nu^\vee)=-1\}$ and for such
$\gamma$ we have $s_\nu(\gamma)=\gamma+\nu >0$. Hence 
$N(w^{-1}_\nu)\cap N(s_\nu)=\varnothing$. The second equality follows from
(i) and the fact that $\ell(w_\nu)=(\rho,\theta^\vee-\nu^\vee)$, see
\cite[Theorem\,4.2]{lp}.

(iii) \ Obvious.
\end{proof}%
In what follows, we write $N(s_\nu)^0$ 
for $\{\gamma\in\Delta^+\mid (\gamma,\nu^\vee)=1 \ \ 
\& \ \ \gamma \prec \nu\}$. 

\begin{s}{Proposition}  \label{char-dom1}  
Suppose $\nu\in\Delta^+_l$. Then
\begin{itemize}
\item[\sf (i)] \ $w_\nu s_0$ is a minimal element with rootlet $\nu$,
and $I_{w_\nu s_0}$ is Abelian. The element $w_\nu s_0$ is also maximal, 
if $\nu\not\in\Pi$.
\item[\sf (ii)]  $s_\nu w_\nu s_0$ is dominant and 
${\rt}(s_\nu w_\nu s_0)=-\nu$. Next, $s_\nu w_\nu s_0\in \HW_{min}$ if 
and only if $\nu\not\in\Pi$;
$s_\nu w_\nu s_0\in \HW_{max}$ if and only if $(\theta,\nu)=0$.
\item[\sf (iii)] \ If $\nu\not\in\Pi$, then
the ideal $I_{s_\nu w_\nu s_0}$ is not Abelian.
\item[\sf (iv)] \ If $\nu\in\Pi$, then
$I_{s_\nu w_\nu s_0}=I_{w_\nu s_0}$ is Abelian.
\end{itemize}
\end{s}\begin{proof}
(i) \ Since $w_\nu s_0(\ap_0)+\delta=\nu$, we obtain ${\rt}(w_\nu s_0)=
\nu$.  
That $w_\nu s_0$ is minimal (in particular, dominant) and the
ideal $I_{w_\nu s_0}$ is Abelian is shown in \cite[Theorem\,4.2]{lp}.
The assertion on maximality stems from 
Proposition~\ref{heis0}.

(ii) \ Set $w=s_\nu w_\nu s_0$. From the very definition
of $w_\nu$, it follows that $w(\ap_0)+\delta=-\nu$. Hence 
the rootlet is as required.
\\[.5ex]
In order to prove that $w$ is dominant, we apply Proposition~\ref{heis0}(i).
Here $v=s_\nu w_\nu$.
If $\ap\in\Pi$ and $(\theta,\ap)=0$, then
 $s_\nu w_\nu(\ap)=w_\nu(\ap)-
(w_\nu(\ap),\nu^\vee)\nu=w_\nu(\ap)>0$.
The last inequality follows from the fact that, by part (i), $w_\nu s_0$ 
is dominant. It is also not hard to give a direct argument.
(Indeed, if $w_\nu(\ap)=-\mu<0$, then $\mu\in N(w_\nu^{-1})$. Therefore
$(\mu,\nu^\vee)=-1$ by Eq.~\re{-1}. That is,
$1=(w_\nu(\ap),\nu^\vee)=(\ap,\theta^\vee)$, a contradiction!
\\[.5ex]
The assertions on maximality and minimality follows from 
Proposition~\ref{heis0}(iii),(iv).

(iii) \ 
Notice that $s_\nu w_\nu s_0(2\delta-\theta)=-\nu <0$. 
Since $s_\nu w_\nu s_0$ is minimal for $\nu\not\in\Pi$, it follows from
Proposition~\ref{bij-mi} that $\theta\in I^2$.

(iv) \ By Eq.~\re{uv} and Lemma~\ref{dlina}(ii), we have 
$N(s_\nu w_\nu s_0)=N(w_\nu s_0)\cup s_0 w_\nu^{-1}(N(s_\nu))$. If $\nu\in
\Pi$, then $N(s_\nu)=\{\nu\}$, and the second set equals 
$\{2\delta-\theta\}$. Thus, $N(s_\nu w_\nu s_0)$ and $N(w_\nu s_0)$ have
the same first layers. Since $I_{w_\nu s_0}$ is Abelian and
$N(w_\nu s_0)=\delta-I_{w_\nu s_0}$, we are done.
\end{proof}%
The next assertion yields a converse to Proposition~\ref{char-dom1},
and completes our characterization of the dominant elements of Heisenberg 
type.

\begin{s}{Proposition}  \label{char-dom2}
Suppose $w=vs_0\in \HW_{dom}$, where $v\in W$. Set
$\nu{=}\left\{\begin{array}{rc} v(\theta), & \text{if \ } v(\theta)>0 \\
-v(\theta), & \text{if \ } v(\theta)<0 \end{array}\right.$.
Then $v=w_\nu$, if $v(\theta)>0$; and
$v=s_\nu w_\nu$, if $v(\theta)=-\nu<0$.
\end{s}\begin{proof}
1. Suppose $\nu=v(\theta)>0$.  
Then $(\gamma,\nu^\vee)=(v^{-1}(\gamma),\theta^\vee)$. 
Therefore, if $(\gamma,\nu^\vee)=-1$, then $v^{-1}(\gamma)<0$.
In view of Eq.~\re{-1}, this means that $N(v^{-1})\supset N(w_\nu^{-1})$.
It follows that there is a "reduced decomposition"
$v^{-1}=\kappa w_\nu^{-1}$ (i.e., $\ell(v)=\ell(\kappa)+\ell(w_\nu)$).
Then $v=w_\nu\kappa^{-1}$ and $\kappa^{-1}(\theta)=\theta$.
This means that $\kappa^{-1}$ takes the (possibly reducible)
root system $\Delta\cap \langle\theta\rangle^\perp$ to itself.
Therefore, if $\kappa\ne \textrm{id}$, then there is an $\ap\in\Pi\cap
\langle\theta\rangle^\perp$
such that $\kappa^{-1}(\ap)<0$. Since 
$v=w_\nu\kappa^{-1}$ is a reduced decomposition and therefore
$N(v)\supset N(\kappa^{-1})$, we have 
$v(\ap)<0$ as well. But Proposition~\ref{heis0}(i) says that $v(\ap)$ must
be positive here.
This contradiction shows that $\kappa=\textrm{id}$,
and we are done.

2.  Suppose $-\nu=v(\theta)<0$.  
We use the same idea as in part 1. If 
$N(v^{-1})\supset N(w_\nu^{-1}s_\nu)$, then we will obtain a reduced
decomposition of the form $v^{-1}=\kappa w_\nu^{-1}s_\nu$, with
$\kappa^{-1}(\theta)=\theta$, and eventually
prove that $\kappa=\textrm{id}$.
So, it suffices to establish the above containment.
By Lemma~\ref{dlina}(ii) and Eq.\re{uv}, we have 
\[
N(w_\nu^{-1}s_\nu)=N(s_\nu)\cup s_\nu{\cdot}N(w_\nu^{-1})=
\{\nu\}\cup \{\gamma\in\Delta^+\mid (\gamma,\nu^\vee)=1 \ \ 
\& \ \ \gamma \prec \nu\} \cup s_\nu{\cdot}N(w_\nu^{-1}) \ .
\]
Let us check that the three sets in the RHS are contained in
$N(v^{-1})$.

$\bullet$ \ $v^{-1}(\nu)=-\theta <0$;

$\bullet$ \ If $(\gamma,\nu^\vee)=1$, then
$1=(v^{-1}(\gamma),-\theta^\vee)$. Hence $v^{-1}(\gamma)<0$;

$\bullet$ \ If $\gamma\in N(w_\nu^{-1})$, then 
$v^{-1}s_\nu(\gamma),\theta^\vee)=(s_\nu(\gamma, -\nu^\vee)=
(\gamma,\nu^\vee)=-1$. Hence $v^{-1}s_\nu(\gamma)<0$.
\end{proof}%
Combining Propositions~\ref{char-dom1} and \ref{char-dom2},
we obtain

\begin{s}{Theorem}   \label{heis-main}
Suppose $w=vs_0$, where $v\in W$. Then $w\in \HW_{dom}$ if and only if
$v$ is either $w_\nu$ or $s_\nu w_\nu$ for some $\nu\in\Delta^+_l$.
Furthermore, ${\rt}(w)$ is $\nu$ (resp. $-\nu$) in the first
(resp. second) case.
\end{s}%
It follows that the number of dominant elements of Heisenberg type is equal
to $\#(\Delta_l)$. Using Proposition~\ref{char-dom1}(ii), we also conclude
that the number of minimal elements of Heisenberg type is equal
to $\#(\Delta_l\setminus\Pi)$, and
the number of maximal elements of Heisenberg type is equal
to $\#(\Delta^+_l\setminus\Pi)+ \#( \{ 
\gamma\in \Delta^+_l\mid (\gamma,\theta)=0\})$. However, it is not yet clear
that if $I$ is a non-trivial ideal (resp. strictly positive ideal) in 
$\gH$, then the corresponding minimal (resp. maximal) element is 
necessarily of Heisenberg type. Actually, this appears to be true for minimal
elements, but not for maximal.

\begin{s}{Proposition}  \label{min-heis}
If $w\in\HW_{min}$ and $I_w\subset \gH$, then 
$w$ is of Heisenberg type.
\end{s}\begin{proof}
It is straightforward to verify that $s_\theta s_0\in \HW_{min}$ 
and $I_{s_\theta s_0}=\gH$ (cf. \cite[Example\,2.6]{duality}).
It follows that $N(w)\subset N(s_\theta s_0)$. Therefore
$s_\theta s_0=\tilde w w$ is a reduced decomposition. Since $w=w's_0$ for 
some $w'\in\HW$, we conclude that $s_\theta=\tilde w w'$ is a reduced 
decomposition for $s_\theta$. Whence $\tilde w, w'\in W$.
\end{proof}%
{\sf Remark.} It happens quite often that $w\in\HW_{max}$ and $I_w\subset
\gH$, but $w$ is not of Heisenberg type. A "uniform" example can be
described as follows. Suppose $\theta$ is a fundamental weight and
$\nu$ is the unique simple root such that $(\theta,\nu)>0$.
According to Proposition~\ref{char-dom1}(ii), $s_\nu w_\nu s_0$ is
dominant but neither maximal nor minimal. Letting $I=I_{s_\nu w_\nu s_0}$,
one can show that $w_{min}(I)=w_\nu s_0$ (this follows from
Proposition~\ref{char-dom1}(iv)) and $w_{max}(I)=s_0 s_\nu w_\nu s_0$
(a straightforward computation).
That is, $w_{max}(I)$ is not of Heisenberg type.
\\[.6ex]
Thus, we have obtained a bijection between 

$\bullet$\quad the roots in $\Delta^+_l \setminus (-\Pi)$;

$\bullet$\quad the minimal elements of Heisenberg type;

$\bullet$\quad 
the non-trivial ideals  lying in $\gH$.
\\
Under this bijection, the roots in $\Delta^+_l\setminus (-\Pi)$ are
obtained 
as the rootlets of the corresponding ideals. In particular, 
different ideals of Heisenberg type have different rootlets.
\\[.5ex]
An explicit construction of the Heisenberg-type elements is given in 
Proposition~\ref{char-dom1}. It is also possible to explicitly
describe the corresponding ideals. 

\begin{s}{Theorem}  \label{main3}  Suppose $\nu\in\Delta^+_l$. Then
\begin{itemize}
\item[\sf (i)] \  $I_{ w_\nu s_0}=\{ \theta\} \cup
\{\theta+ w^{-1}_\nu(N(w^{-1}_\nu) \}= \{\theta\} \cup
\{ \theta-N(w_\nu) \}$;
\item[\sf (ii)] \ $I_{s_\nu w_\nu s_0}=\{\theta\} \cup
\{ \theta-N(w_\nu)\}\cup \{ \theta- w^{-1}_\nu(N(s_\nu)^0) \}=
I_{ w_\nu s_0}\cup \{ \theta- w^{-1}_\nu(N(s_\nu)^0) \}$, if $\nu\not\in
\Pi$.
\end{itemize}
\end{s}\begin{proof}
(i) Since $I_{ w_\nu s_0}$ is Abelian, $\#(I_{ w_\nu s_0})=
\ell(w_\nu s_0)=(\rho,\theta^\vee-\nu^\vee)+1$.
Since the set in the right-hand side has the required cardinality,  
it suffices to verify that it is contained in
$I_{ w_\nu s_0}$. Observe that
$(\gamma,\theta^\vee)=(w_\nu(\gamma),\nu^\vee)=-1$ for any
$\gamma\in w^{-1}_\nu(N(w^{-1}_\nu)$. Therefore we know how 
$s_0$ acts on $w^{-1}_\nu(N(w^{-1}_\nu))$. Then
\[
w_\nu s_0(\delta-\theta-w^{-1}_\nu(N(w^{-1}_\nu))=
w_\nu(-w^{-1}_\nu(N(w^{-1}_\nu))=-N(w^{-1}_\nu) \subset -\Delta^+ \ .
\]
Also, $w_\nu s_0(\delta-\theta)=\nu-\delta <0$.

(ii) Remember that $\gH^2=
\{\theta\}$. Since $I:=I_{s_\nu w_\nu s_0}\subset \gH$ and is not
Abelian by Proposition~\ref{char-dom1}(iii), 
we have $I^2=\{\theta\}$ as well. Therefore, using 
Corollary~\ref{ell}, we obtain
$\#(I)=\ell(s_\nu w_\nu s_0)-1=\ell(s_\nu w_\nu)=
(\rho,\theta^\vee+\nu^\vee)-1$.
Here, again, the set in the right-hand side has the prescribed cardinality, 
Hence it remains to show that it is contained in $I$. 
First, it follows from \re{uv} that $I_{w_\nu s_0}\subset
I_{s_\nu w_\nu s_0}$. Next, we have 
$(w^{-1}_\nu(N(s_\nu)^0),\theta^\vee)=(N(s_\nu)^0, \nu^\vee)=1$.
Therefore, we know how $s_0$ acts on $w^{-1}_\nu(N(s_\nu)^0)$, and hence
we can compute that all roots in
$s_\nu w_\nu s_0(\delta-\theta+w^{-1}_\nu(N(s_\nu)^0))$ are negative.
\end{proof}%
{\sf Remark.} In principle, it is harmless to omit the hypothesis in 
part (ii) that $\nu\not\in\Pi$. For, $N(s_\nu)^0=\varnothing$, if $\nu\in
\Pi$. Then part (ii) would assert in this case
that $I_{s_\nu w_\nu s_0}=I_{w_\nu s_0}$,
which was already proved in Proposition~\ref{char-dom1}(iv).


\section{Minimax elements and ideals: general properties}  
\label{minm}
\setcounter{equation}{0}

\noindent
An element $w\in \HW$ is called {\it minimax\/}, if it is simultaneously
minimal and maximal. Combining the corresponding definitions of
Section~\ref{old}, we obtain the following:

\begin{rem}{Definition}  \label{def-mm}
$w\in\HW$ is called {\it minimax\/}, if
\begin{itemize}
\item[\sf (i)] \ $w$ is dominant;
\item[\sf (ii)] \ if $\ap\in\HP$ and  
$w^{-1}(\ap)=k\delta+\mu$ for some $\mu\in \Delta$, then \ $-1\le k\le 1$.
\end{itemize}
\end{rem}%
The corresponding (strictly  positive) ideal is said to be
{\it minimax\/}, too. That is, $I\in \AD_0$ is minimax, 
if $w_{min}(I)=w_{max}(I)$. Therefore we merely write $w(I)$, if 
$I$ is minimax.
The set of all minimax elements, which is clearly equal
to $\HW_{min}\cap\HW_{max}$, is denoted by 
$\HW_{mm}$. Accordingly, $\AD_{mm}$ is the set of minimax ideals
in $\Delta^+$.
\\[.6ex]
Combining geometric descriptions of minimal and maximal elements given 
in Section~\ref{old}, we obtain the following. 
\begin{s}{Proposition} 
\label{opis-mm} \\
Set $ D_{mm}=D_{min}\cap D_{max}=\{
x\in V\mid -1\le (x,\ap)\le 1 \ \ \forall\ap\in\Pi \ \ \& \ \ 
0\le (x,\theta)\le 2\}$.  Then \par
1. $w=v{\cdot}t_r\in \HW_{mm}$ \ $\Longleftrightarrow
\left\{ \begin{array}{l} w \text{ is dominant},  \\
                         v(r)\in D_{mm}\cap {Q}^\vee \ .
       \end{array}\right.$

2. The mapping $\HW_{mm} \to D_{mm}\cap Q^\vee$, 
$w=v{\cdot}t_r\mapsto v(r)$, is a bijection.
\end{s}%
It can be shown that there is an element of $\HW$ that
takes $D_{mm}$ to 
\[
 \{x\in V\mid 0\le (x,\ap)\le 2 \ \ \forall\ap\in\Pi \ \ \& \ \ 
h-1\le (x,\theta)\le h+1\} \ .
\] 
But, I do not see how it could help us to count the integral points
in $D_{mm}$.
\\
Recall that, for any strictly positive ideal $I\subset\Delta^+$ and 
any $\gamma\in I$, we have defined the numbers 
\[
 l(\gamma,I)=\max\{ m\mid \gamma=\sum_{i=1}^m \varkappa_i,\ \ 
   \varkappa_i\in I\}\ ,\quad
 k(\gamma,I)=\min\{ n\mid \gamma=\sum_{i=1}^n \nu_i,\ \ 
   \nu_i\in \Delta^+\setminus I\} \ .
\]
Using results of Section~\ref{old} we immediately obtain various 
reformulations of the minimax property. 

\begin{s}{Theorem}  \label{pereform} 
Suppose $I$ is a strictly positive ideal and $w$ is 
either $w_{min}(I)$ or $w_{max}(I)$.
Then the following conditions are equivalent:
\begin{itemize}
\item[\sf (i)]  $w\in\HW_{mm}$ ;
\item[\sf (ii)]  $k(\gamma,I)-1=l(\gamma,I)$ for all $\gamma\in I$;
\item[\sf (iii)] $R_I=w^{-1}\asst\gA$;
\item[\sf (iv)]  $R_I$ consists of a single alcove.
\end{itemize}
\end{s}\begin{proof}
Use Eq.~\re{l<k} and Theorem~\ref{min+max}.
\end{proof}%
Thus, the minimax elements (ideals)
are in a bijection with the dominant regions of the
Catalan arrangement consisting of a single alcove.
\\[.5ex]
Now, we give a description of
minimax Abelian ideals.
Let $I\subset \Delta^+$ be a nontrivial Abelian ideal
and $w=w_{min}(I)$ the corresponding minimal element. It was shown in
\cite[Prop.\,2.5]{lp} that $w(\ap_0)=-\delta+\nu$, where $\nu\in\Delta^+_l$.
In particular, ${\rt}(w)={\rt}(I)$ is a positive root. 

\begin{s}{Theorem}  \label{abel-mm}
An Abelian ideal $I$ is minimax if and only if\/ ${\rt}(I)$
is not a simple root.
\end{s}\begin{proof}
1. Assume $w(\ap_0)+\delta=\ap_i\in\Pi$. Then 
$w^{-1}(\ap_i)=2\delta-\theta$, i.e., $w$ is not maximal.

2. Since $I$ is Abelian, $l(\gamma,I)=1$ for all $\gamma\in I$. 
Therefore, in view of Theorem~\ref{pereform}(ii), the minimax property
is equivalent to that any $\gamma\in I$ is a sum of two elements of
$\Delta^+\setminus I$.
Assume $w(\ap_0)+\delta=\mu\not\in\Pi$. Then
$\mu=\mu_1+\mu_2$, where $\mu_i\in\Delta^+$. Hence $w^{-1}(\mu_i)\in
\HD^+$ and $2\delta-\theta=w^{-1}(\mu_1)+w^{-1}(\mu_2)$ is a sum of
two positive roots. The only possibility for this is 
$\theta=\gamma_1+\gamma_2$ and
$w^{-1}(\mu_i)=\delta-\gamma_i$. Since $w(\delta-\gamma_i)=\mu_i \succ 0$,
we have $\gamma_i\not\in I$. Thus, the required property is satisfied 
for $\theta\in I$. To prove this for any $\mu\in I$, we use a
descending induction. Indeed, if $\theta-\ap\in I$ for some
$\ap\in\Pi$, then the condition 
$\gamma_1+\gamma_2-\ap \in I\subset \Delta^+$ implies that 
$\ap\ne\gamma_i$ and at least one of
$\gamma_i-\ap$ is in $\Delta^+$ (and hence in $\Delta^+\setminus I$),
see \cite[Lemma\,2.3]{duality}. Hence $\theta-\ap$ is a sum of two
elements from $\Delta^+\setminus I$. This can be continued further.
\end{proof}%
It should already be clear that a minimax ideal is not necessarily
Abelian. The simplest example is the ideal in ${\frak sl}_5$
generated by $\{\ap_1+\ap_2,\, \ap_3+\ap_4\}$.
More precisely,
applying results of Section~\ref{heis} to minimax elements, we obtain

\begin{s}{Proposition}  \label{heis-mm}
A Heisenberg type element $w=vs_0$ is minimax if and only 
if either $v=w_\nu$, where $\nu\in\Delta^+_l\setminus\Pi$, or
$v=s_\nu w_\nu$, where $\nu\in\Delta^+_l\setminus\Pi$ and 
$(\theta,\nu)=0$.
The minimax ideals of the second kind are not Abelian.
\end{s}%
It is natural to look at $\Delta_{mm}$, the set of rootlets of
all minimax elements (ideals). Clearly, we have the inclusion
$\{ \nu\in \Delta_l \mid (\nu,\theta)\ge 0\}\setminus \{\pm\Pi\}
\subseteq \Delta_{mm}$. But, in general, it is a strict containment.
This can be expressed in the following way.
Consider the mapping 
$\zeta: \AD_{mm}\to \Delta_{mm}$ that associates the rootlet
to a minimax ideal. Then, in general, there exist fibres of $\zeta$
that do not contain a Heisenberg type ideal.
Here is an example. Let $\g={\frak sl}_6$. Take the ideal generated
by the roots $\ap_1{+}\ap_2,\,\ap_3{+}\ap_4$. A straightforward computation
shows that it is minimax and the rootlet is $\nu=-(\ap_1{+}\ap_2{+}\ap_3)$.
Thus, $\nu$ is negative and $(\nu,\theta)\ne 0$.

\noindent
It is intuitively clear that if a (dominant) region of $\Cat(\Delta)$
consists of a single alcove, then this region must be sufficiently
close to the origin. Equivalently, the corresponding ideal must be
not too large. We now discuss a precise meaning that can be given
to this intuitive feeling. 
Given $I\in \AD$, consider the quantity $\#\Gamma(I)+\#\Xi(I)$.
It is not hard to realize that its maximal value  is
$2p-1$, which is attained for $I=\Delta^+\setminus \Pi$.

\begin{s}{Proposition}  \label{bound-mm}
Suppose $I\in\AD_{mm}$; then 
\begin{itemize}
\item[\sf (i)] \ $\#\Gamma(I)+\#\Xi(I)\le p+1$;
\item[\sf (ii)] \ If \ $\#\Gamma(I)+\#\Xi(I)= p+1$, then $h$ is odd.
\end{itemize}
\end{s}\begin{proof}
(i) \  
Set $w=w(I)$ and write  $w^{-1}(\ap)=k_\ap\ap +\nu_\ap$.
Since $w$ is simultaneously minimal and maximal, $k_\ap\in\{-1,0,1\}$ for
any $\ap\in\HP$. (Notice that $\nu_\ap\in \Delta^+$, if $k_\ap=-1$; 
$\nu_\ap\in -\Delta^+$, if $k_\ap=1$.)
Making use of Propositions~\ref{gen-mi} and \ref{gen-ma}
we obtain

$\gamma\in \Gamma(I) \ \Longleftrightarrow \ \text{there is an }\ \ap\in\HP
\text{ such that }\ w^{-1}(\ap)=-\delta+\gamma$,

$\mu\in \Xi(I)\ \Longleftrightarrow \ \text{there is an }\ \ap\in\HP
\text{ such that }\ w^{-1}(\ap)=\delta-\mu$.

\noindent
This means that the elements of $\Gamma(I)$ (resp. $\Xi(I)$)
are in a bijection with the affine simple roots such that 
$k_\ap=-1$ (resp. $k_\ap=1$). Since $\#\HP=p+1$, we are done.

(ii) \  It follows from the previous
part of proof that  if $\#\Gamma(I)+\#\Xi(I)= p+1$, then $k_\ap\ne 0$
for all $\ap$. Therefore we obtain the partition $\HP=\HP_+\sqcup\HP_-$,
according to whether $k_\ap=+1$ or $-1$. 
Recall that the coefficients $c_\ap$ $(\ap\in\Pi)$ are defined by the
equality $\theta=\sum_{\ap\in\Pi}c_\ap\ap$. We also set $c_\ap=1$ for 
$\ap=\ap_0$. Then
\[
   \delta=\sum_{\ap\in \HP}c_\ap\ap  \quad\text{ and } \quad \
   \sum_{\ap\in\HP}c_\ap=h \ .
\]
Because also $\delta=\sum c_\ap w^{-1}(\ap)$, we obtain
\[
  1=\sum_{\ap\in\HP_+}c_\ap-\sum_{\ap\in\HP_-}c_\ap=h-
  2\sum_{\ap\in\HP_-}c_\ap \ .
\]
That is, $h$ is odd.
\end{proof}%
As is well-known, $h$ is odd if and only if $\g={\frak sl}_{2n+1}$.
The following example shows that the above upper bound is exact.
Consider the ideal in ${\frak sl}_{2n+1}$ such that
\[
  \Gamma(I)=\{\ap_1{+}\ap_2,\ap_3{+}\ap_4,\dots,\ap_{2n-1}{+}\ap_{2n}\} \ .
\]
Then 
\[
\Xi(I)=\{\ap_1,\ap_2{+}\ap_3,\dots,\ap_{2n-2}{+}\ap_{2n-1},\ap_{2n}\} \ .
\]
It follows from the description of minimax ideals in type $\GR{A}{p}$, see
Corollary~\ref{nm-mm-sl}, that this ideal is minimax.
For all other simple Lie algebras, it is possible to exhibit a minimax
ideal such that $\#\Gamma(I)+\#\Xi(I)= p$.


\section{Counting minimax elements/ideals in the simple Lie algebras}  
\label{count}
\setcounter{equation}{0}

\noindent
By Proposition~\ref{opis-mm}, the number of minimax elements of
$\HW$ is equal to the cardinality of
$D_{mm}\cap Q^\vee$. Unfortunately, $D_{mm}$ is not a simplex, 
so that it is not easy to find a uniform expression for
$\#(D_{mm}\cap Q^\vee)$. However, for each simple algebra one has a
certain system of inequalities and one may try to solve these systems
individually.
It turns out that, for practical computations, it is easier to deal with
the coweight lattice in $V$, denoted $P^\vee$. The number
$[P^\vee: Q^\vee]$ is called the {\it index of connection\/} of $\Delta$.

%
\noindent
Recall from Section~\ref{old}
that $\{\varpi_i^\vee\}_{i=1}^p$ is the dual basis of $V$ for 
$\{\ap_i\}_{i=1}^p$. Then the lattice generated by the
$\varpi_i$'s is $P^\vee$. 
If $y=\sum_{i} y_i\varpi_i\in P^\vee$, then $y\in Q^\vee$ if and only if
a certain congruence condition is satisfied for $(y_1,\dots,y_p)\in
{\Bbb Z}^p$. 
Recall that $\theta=\sum_{i=1}^p c_i\ap_i$. 
The equations of $D_{mm}$ in Proposition~\ref{opis-mm}
show that $\#(D_{mm}\cap P^\vee)$ equals the number of solutions of
the following system:
\begin{equation}  \label{syst}
    \left\{ \begin{array}{c}
          y_i\in\{-1,0,1\} \ \ (i=1,2,\dots,p) \\
          0\le c_1y_1+\ldots +c_py_p\le 2  
\ \  .
\end{array}\right.
\end{equation}
It is convenient to introduce one more variable, $y_0$,
which also ranges over $\{-1,0,1\}$, and add it to the 
expression in the second inequality, so that the total sum to
be equal to 1. (We also set $c_0=1$.) Then the above system becomes
\begin{equation}  \label{syst-ext}
    \left\{ \begin{array}{c}
          y_i\in\{-1,0,1\} \ \ (i=0,1,\dots,p) \\
          c_0y_0+c_1y_1+\ldots +c_py_p=1  
\ \  .
\end{array}\right.
\end{equation}
In a sense, this procedure corresponds to taking the extended 
Dynkin diagram of $\g$. For this reason, system~\re{syst-ext} will
be referred to as the {\it extended\/} system.
\\
It is easily seen that the number of solutions of system \re{syst-ext}
is nothing but the coefficient of $x$ in the expansion of the Laurent
polynomial 
\begin{equation}  \label{laurent}
\prod_{i=0}^p(x^{-c_i}+1+x^{c_i}) \ .
\end{equation}
An important particular case is that where all $c_i=1$. 
The coefficients in the expansion of $(x^{-1}+1+x)^n$ are called
{\it trinomials\/}, and we write $X_k(n)$ for the coefficient of
$x^k$. The coefficient of $x^0$ (resp. $x$) is called the {\it central\/}
(resp. {\it next-to-central\/})
trinomial. It is easily  seen that
\[
  X_k(n):=\sum_{l\ge 0}\frac{ n!}{l! (k+l)! (n-2l-k)!} \ . 
\]
The trinomial coefficients satisfy the following relation
\begin{equation}  \label{tri-ind}
   X_k(n{+}1)=X_{k-1}(n)+X_k(n)+x_{k+1}(n)  \ .
\end{equation}
To get the number $\#(D_{mm}\cap Q^\vee)$,
one should add to both systems the above-mentioned congruence condition.
Recall that $[P^\vee: Q^\vee]$ is equal to 
the number of 1's among the coefficients $c_i$, $i=0,1,\dots,p$.
Our first goal is to show that the congruence condition can be omitted.

\begin{s}{Theorem}   \label{P:Q}
$\#(D_{mm}\cap P^\vee)=
[P^\vee: Q^\vee]\cdot\#(D_{mm}\cap Q^\vee)$.
\end{s}\begin{proof}
The only proof I know is case-by case (see also Question~4 in
Section~\ref{concl}). In considering particular cases, our numbering of
$\ap_i$'s and hence of $c_i$'s correspond to the numbering adopted in
\cite[Tables]{VO}. The explicit form of the congruence condition can be extracted
from Table~3 in loc.\,cit.

{\bf 1)} \ $\g={\frak sl}_{n+1}$. \quad
Here all $c_i=1$ and $[P^\vee:Q^\vee]=n+1$. The extended system is
\begin{equation}  \label{syst-sl}
\left\{ \begin{array}{c}
y_0+y_1+\ldots+y_n=1 \\ 
y_i\in\{-1,0,1\} \ , \end{array}\right.  
\end{equation}
and the congruence is \ 
$ny_1+(n-1)y_2+\ldots +y_n \in (n+1){\Bbb Z}$. 
Obviously, the symmetric group ${\Bbb S}_{n+1}$ acts on the set of 
solutions
of system~\re{syst-sl}. Let $\tilde c\in {\Bbb S}_{n+1}$ be the
transformation $y_i\mapsto y_{i-1}$ ($i>0$),
$y_0\mapsto y_n$. (It is a Coxeter element of ${\Bbb S}_{n+1}$.)
It is easily seen that all orbits of $\tilde c$ on the set of solutions of
Eq.~\re{syst-sl} have cardinality 
$n+1$. Furthermore, making use of the relation $y_0+y_1+\ldots+y_n=1$,
one easily verifies that
\[
\tilde c(ny_1+(n-1)y_2+\ldots +y_n)- (ny_1+(n-1)y_2+\ldots +y_n)
\equiv 1 \pmod{n+1} \ . 
\]
This yields the desired relation.

{\bf 2)} \ $\g={\frak sp}_{2n}$. \quad
Here $[P^\vee:Q^\vee]=2$, the extended system is
\begin{equation}  \label{syst-sp}
\left\{ \begin{array}{c}
y_0+2(y_1+\ldots+y_{n-1})+y_n=1 \\ 
y_i\in\{-1,0,1\} \ , \end{array}\right.  
\end{equation}
and the congruence is \ 
$y_n\in 2{\Bbb Z}$. It follows that $y_0+y_n$ is always odd 
for a solution of \re{syst-sp} and
that the permutation of $y_0$ and $y_n$ takes the solutions  
with $y_n$ even to those with $y_n$ odd. 
Hence the result.

{\bf 3)} \ $\g={\frak so}_{2n+1}$. \quad
Here $[P^\vee:Q^\vee]=2$, the extended system is
\begin{equation}  \label{syst-so}
\left\{ \begin{array}{c}
y_0+y_1+2(y_2+\ldots+y_{n})=1 \\ 
y_i\in\{-1,0,1\} \ , \end{array}\right.  
\end{equation}
and the congruence is \ 
$y_1+y_3+y_5+\ldots \in 2{\Bbb Z}$. Here one may use the permutation
of $y_0$ and $y_1$, as in part 2).

{\bf 4)} \ $\g={\frak so}_{2n}$, $n\ge 4$. \quad
Now $[P^\vee:Q^\vee]=4$ and the extended system is
\begin{equation}  \label{syst-so-even}
\left\{ \begin{array}{c}
y_0+y_1+2(y_2+\ldots+y_{n-2})+y_{n-1}+y_n=1 \\ 
y_i\in\{-1,0,1\} \ . \end{array}\right.  
\end{equation}
But the congruence condition depends on the parity of $n$.

(A) $n=2l$. \ Here we actually have two conditions:
\begin{equation}  \label{congr-even}
   \left\{ \begin{array}{r}  y_{n-1}+y_n \in 2{\Bbb Z}  \\
         y_1+y_3+\ldots +y_{n-1} \in 2{\Bbb Z} \end{array}\right.  \ .
\end{equation}
Here we consider the cyclic permutation $\tilde c:\ y_0\to y_1\to y_{n-1}\to
y_n\to y_0$. Since the sum $y_0+y_1+y_{n-1}+y_n$ is always odd, all
orbits of this permutation on the set of solutions of 
system~\re{syst-so-even} have cardinality 4. It is not hard to verify
that each orbit contains a unique representative satisfying 
\re{congr-even}. Let us give some details.
It follows from system~\re{syst-so-even} that \ 
$-1\le y_2+\ldots +y_{n-2}\le 2$. That is, there are four
cases:

$(a_1)$ \ \ $y_2+\ldots +y_{n-2}=-1$  and hence $y_0+y_1+y_{n-1}+y_n=3$;

$(a_2)$ \ \ $y_2+\ldots +y_{n-2}=0$ \ and hence $y_0+y_1+y_{n-1}+y_n=1$;

$(a_3)$ \ \ $y_2+\ldots +y_{n-2}=1$ \ and hence $y_0+y_1+y_{n-1}+y_n=-1$;

$(a_4)$ \ \ $y_2+\ldots +y_{n-2}=2$ \ and hence $y_0+y_1+y_{n-1}+y_n=-3$;
\\[.5ex]
The number of possibilities for $(y_0,y_1,y_2,y_3)$ in these four cases is 
equal to 4,\,16,\,16, and 4, respectively. 
(For instance, in case $(a_2)$, this number is $X_1(4)=16$.)
Hence, the total number is 40, and one have to test only 10 orbits of 
$\tilde c$ (actually, only five orbits in view of the symmetry).
Each orbit contains two representatives satisfying the first condition
in \re{congr-even}; for these two representatives, the parity of
$y_1$ is different, so that only one of them satisfies
the second condition.
Notice also that we have proved that
the total number of the solutions of system~\re{syst-so-even} is equal to
\[
    4X_{-1}(n-3)+16X_{0}(n-3)+16X_{1}(n-3)+4X_{2}(n-3) \ .
\]

(B) $n=2l+1$. \ Here we have one congruence condition:
\[  
   2(y_1+y_3+\ldots +y_{n-2})+y_{n-1}-y_n \in 4{\Bbb Z} \ .
\]  
However, we can split it in two conditions:
$y_{n-1}-y_n \in 2{\Bbb Z}$  and 
$y_1+y_3+\ldots +y_{n-2}+\frac{1}{2}(y_{n-1}-y_n) \in 2{\Bbb Z}$, so that
the previous argument goes through verbatim.

{\bf 5)} \ $\g=\GR{E}{6}$. \quad
Here $[P^\vee:Q^\vee]=3$, the extended system is
\begin{equation}  \label{syst-e6}
\left\{ \begin{array}{c}
y_0+y_1+2y_2+3y_3+2y_4+y_5+2y_{6}=1 \\ 
y_i\in\{-1,0,1\} \ , \end{array}\right.  
\end{equation}
and the congruence condition is \ 
$y_1-y_2+y_4-y_5 \in 3{\Bbb Z}$. Then one may use the permutation
\[
\tilde c:\ \left\{ \begin{array}{c}
y_1\to y_5 \to y_0\to y_1 \\
y_2\to y_4 \to y_6 \to y_2 
 \end{array}\right.  \ ,
\]
since $y_1-y_2+y_4-y_5 -\tilde c(y_1-y_2+y_4-y_5)
 \equiv 1 \pmod{3}$

{\bf 6)} \ $\g=\GR{E}{7}$. \quad
Here $[P^\vee:Q^\vee]=2$, the extended system is
\begin{equation}  \label{syst-e7}
\left\{ \begin{array}{c}
y_0+y_1+2y_2+3y_3+4y_4+3y_5+2y_{6}+2y_7=1 \\ 
y_i\in\{-1,0,1\} \ , \end{array}\right.  
\end{equation}
and the congruence condition is \ 
$y_1+y_3+y_7 \in 2{\Bbb Z}$. Then one may use the involution
\[
\tilde c:\ \left\{ \begin{array}{c}
y_1\to y_0, \ y_2\to y_6 \\
y_3\to y_5, \ y_4 \to y_4, \ y_7\to y_7 \ . 
\end{array}\right.  
\]
since $y_1+y_3+y_7 -\tilde c(y_1+y_3+y_7)$ is odd.
\\[.6ex]
{\bf 7)} \ Finally, we have $P^\vee=Q^\vee$  
for $\GR{G}{2},\GR{F}{4},\GR{E}{8}$, and the theorem is proved.
\end{proof}%

\noindent
Now, we are prepared to compute the number of minimax elements in all
simple Lie algebras.
For the classical Lie algebras, we obtain some well-known
combinatorial quantities. 

\begin{s}{Theorem}   \label{sl}
The number of minimax elements in $\HW({\frak sl}_{n+1})$ is equal to
the $n$-{\sf th Motzkin number\/}, $M_n$.
\end{s}%
\vskip-1ex
\begin{s}{Theorem}   \label{sp}
The number of minimax elements in $\HW({\frak sp}_{2n})$  
or $\HW({\frak so}_{2n+1})$ is equal to the number of\/
{\sf directed animals of size\/} $n$, denoted $\textrm{{\mbox dir}}_n$.
\end{s}%
\vskip-1ex
\begin{s}{Theorem}   \label{so}
The number of minimax elements in $\HW({\frak so}_{2n})$  
is equal to $2\textrm{dir}_{n-2}+\textrm{dir}_{n-1}$.
\end{s}%
We refer to \cite{martin}, \cite{dsh}, and \cite[Ex.\,6.37]{rstan} 
for basic facts on Motzkin numbers and to 
\cite{anim}, \cite[Ex.\,6.46]{rstan} for "directed animal" results;
see also explicit formulae below.
The first few terms of these sequences are

\begin{center}
\begin{tabular}{c|cccrrrrr}
    $n$        & 1 & 2 & 3 &  4 & 5  & 6  &   7 &   8\\ \hline
$M_n$          & 1 & 2 & 4 &  9 & 21 & 51 & 127 & 323\\ 
$\text{dir}_n$ & 1 & 2 & 5 & 13 & 35 & 96 & 267 & 750\\
\end{tabular}
\end{center}

\noindent 
{\it Proof of Theorem~\ref{sl}}. \\
Consider the extended system \re{syst-sl}.
Here the number of solutions is just 
\[
   X_1(n{+}1)=\sum_{k\ge 0}\frac{ (n+1)!}{k! (k+1)! (n-2k)!} \ .
\]
It then follows from Theorem~\ref{P:Q} that the number of minimax 
elements is
\[
   \frac{X_1(n{+}1)}{n+1}=\sum_{k\ge 0}\genfrac{(}{)}{0pt}{}{n}{2k}C_k \ ,
\]
where $C_k=\frac{1}{k+1}\genfrac{(}{)}{0pt}{}{2k}{k}$ is the $k$-th
Catalan number. But the right-hand-side is a well-known formula for
$M_n$, see \cite{dsh} and \cite[Eq.\,(21)]{anim}. \qu

{\it Proof of Theorem~\ref{sp}}. \\
In both cases, the extended systems~\re{syst-sp}, \re{syst-so} are the same,
up to renumbering $y_i$. 
Let us work with \re{syst-sp}.
Here we have the constraint 
$0\le y_1+y_2+\ldots +y_{n-1}\le 1$. That is, there are two cases:

$(a_1)$ \ \ $y_1+y_2+\ldots +y_{n-1}=0$ and hence $y_0+y_n=1$;

$(a_2)$ \ \ $y_1+y_2+\ldots +y_{n-1}=1$ and hence $y_0+y_n=-1$.
\\[.5ex]
This means that system~\re{syst-sp} has $2X_0(n-1)+2X_1(n-1)$ solutions.
Thus, the number of minimax elements is
the sum of central and next-to-central trinomials
\begin{multline*}   
   X_0(n-1)+X_1(n-1)=\sum_{k\ge 0}\frac{(n-1)!}{k! k! (n-2k-1)!} +  
 \sum_{k\ge 0}\frac{(n-1)!}{k!(k+1)! (n-2k-2)!}= \\
 =\sum_{k\ge 0}\genfrac{(}{)}{0pt}{}{2k}{k}\genfrac{(}{)}{0pt}{}{n-1}{2k}
+\sum_{k\ge 0}\genfrac{(}{)}{0pt}{}{2k+1}{k}\genfrac{(}{)}{0pt}{}{n-1}{2k+1}=
  \sum_{q\ge 0}\genfrac{(}{)}{0pt}{}{q}{[q/2]}\genfrac{(}{)}{0pt}{}{n-1}{q} \ ,
\end{multline*}
which is a well-known expression for $\text{dir}_n$,
see \cite[Eq.\,(27)]{anim}.  \qus

Notice that we have also derived the relation
\begin{equation} \label{dir-tri}
 X_0(n-1)+X_1(n-1)=\text{dir}_n  \ .
\end{equation}
{\it Proof of Theorem~\ref{so}}. \\
All essential work is already done in the proof of Theorem~\ref{P:Q},
part 4. Namely, the number of minimax elements is one fourth of the
total number of solutions of system~\re{syst-so-even}, i.e.,
\[
   \frac{1}{4}\Bigl( 4X_{-1}(n-3)+16X_{0}(n-3)+16X_{1}(n-3)+4X_{2}(n-3)
\Bigr) \ .
\]
Making use of Eq.~\re{tri-ind} and \re{dir-tri}, one transforms this 
sum to $2\text{dir}_{n-1}+\text{dir}_{n-2}$.
\qus
{\sf Remarks.}
1.
Although $\#\AD_{mm}({\frak sp}_{2n})=\#\AD_{mm}({\frak so}_{2n+1})$,
the ideals occurring in these two sets have different algebraic properties.
For instance, $\AD_{mm}({\frak sp}_{8})$ contains 8 Abelian ideals (of 13), 
whereas $\AD_{mm}({\frak so}_{9})$ contains 11.

2. Using the relation $\text{dir}_n=3\,\text{dir}_{n-1}-M_{n-2}$, we may 
also write $\#\AD_{mm}({\frak so}_{2n})=5\,\text{dir}_{n-2}-M_{n-3}$.

\begin{s}{Theorem}   \label{except}
1.  For $\g=\GR{E}{p}$, $p=6,7,8$,  the number 
of minimax elements in $\HW$ equals 67, 217, and 834, respectively.

2.  For $\g=\GR{F}{4}$ and $\GR{G}{2}$, there are 17 and 3 minimax elements,
respectively.
\end{s}\begin{proof}
In each case, one can directly compute the coefficient 
of $x$ in the respective Laurent polynomial \re{laurent}, which gives 
the number
$\#(D_{mm}\cap P^\vee)$, and then use Theorem~\ref{P:Q}.
\end{proof}%
\vskip-1.2ex
\begin{rem}{Example}
For $\g=\GR{F}{4}$, the coweight and coroot lattices are equal
and system \re{syst} is
\[   
    \left\{ \begin{array}{c}
          y_i\in\{-1,0,1\} \ \ (i=1,\dots,4)\ ,\\
          0\le 2y_1+4y_{2}+3y_3+2y_4\le 2   \ \  .
\end{array}\right.
\]  
It has 17 solutions, and a description of the corresponding 
non-trivial Abelian minimax ideals can be extracted from \cite[Table\,1]{lp},
using Theorem~\ref{abel-mm} (there are 12 such ideals).
The  non-abelian minimax ideals are described in the next table.

\begin{center}
\begin{tabular}{ccccc}
$\Gamma(I)$  & $\#(I)$ & $\#(I^2)$ &  $w(\ap_0)$ & $y$  \\ \hline
$[1,2,1,1]$ & 9  & 1  & $-\delta-[0,2,1,0]$  & $(1,-1,0,1)$ \\
$[1,1,1,1]$ & 10 & 1  & $-\delta-[2,2,1,0]$ & $(-1,0,0,1)$ \\
$[0,2,2,1],\ [2,2,1,0]$ & 10 & 2 & $-2\delta+[2,4,2,1]$ & $(1,1,-1,-1)$ \\
$[0,2,1,1],\ [2,2,1,0]$ & 12 & 3 & $-2\delta+[2,2,1,0]$ & $(0,1,0,-1)$ \\
\hline
\end{tabular}
\end{center}
For all cases we have $I^3=\varnothing$, so that $\ell(w(I))=\#(I)+\#(I^2)$.
We write $[n_1, n_2, n_3, n_4]$ for the root $\sum_i n_i\ap_i$.
The numbering of simple roots follows \cite[Tables]{VO}. The last column
presents the corresponding solution of the above system. The ideals given
in the first two rows are of Heisenberg type, and the explicit description
of $w(I)$ for them can be extracted from Proposition~\ref{char-dom1}(ii)
and Theorem~\ref{main3}(ii).
\end{rem}%
%

\section{A description of the minimax ideals for ${\frak sl}_{n+1}$
and ${\frak sp}_{2n}$}  
\label{sln}
\setcounter{equation}{0}

\noindent
In this section, we give an explicit matrix description of the minimax ideals
in ${\frak sl}_{n+1}$ and ${\frak sp}_{2n}$.
The idea is as follows. In both cases, we describe a certain set of
ideals containing the minimax ones. Then we count the number of
ideals in this, a priori, larger set. It turns out that we again obtain 
the Motzkin numbers and the number of directed animals of size $n$,
respectively, which solves the problem. 

\un{$\g={\frak sl}_{n+1}$}. \\
We choose $\be$ (resp. $\te$) to be the set of upper-triangular (resp.
diagonal)  matrices of size $n+1$. With the usual numbering of rows and columns
of $n{+}1$ by $n{+}1$ matrices, we identify the positive roots with the pairs
$(i,j)$, where $1\le i<j\le n+1$. Here
$\ap_i=(i,i+1)$ and therefore $(a,b)=\ap_a+\ldots+\ap_{b-1}$.
Then an antichain in 
$\Delta^+({\frak sl}_{n+1})$ is an arbitrary sequence 
$\Gamma=\{(a_1,b_1),\ldots,(a_k,b_k)\}$, where $a_i< b_i$ and 
the two sequences
$\{a_i\}$ and $\{b_i\}$ are strictly increasing. 

\begin{rem}{Definition}   \label{not-meet}
Let us say that the generators of an ideal $I\in\AD({\frak sl}_{n+1})$ 
{\it do not meet\/} or $I$ has {\it non-meeting generators\/},
if the antichain $\Gamma(I)$ has the property that $b_j\ne a_i+1$ for 
all $i,j$.
\end{rem}%
%
The equality $b_j=a_i+1$ means 
that the last root in the expression for $(a_j,b_j)$ is equal to
the first root in the expression for $(a_i,b_i)$. This explains
the term.

\begin{s}{Proposition}   \label{nm-sl}
If $I$ is a minimax ideal, then its generators do not meet.
\end{s}\begin{proof*}
Suppose $I$ has meeting generators; i.e., 
$\Gamma(I)\supset \gamma_1,\gamma_2$, where
$\gamma_1=\ap_i+\ap_{i+1}+\ldots +\ap_j$ and
$\gamma_2=\ap_j+\ap_{j+1}+\ldots +\ap_k$
($i<j<k$).
Consider $\gamma:=\gamma_1+\gamma_2-\ap_j=
\ap_i+\ldots +\ap_k$. It is easily seen that, for any presentation
of $\gamma$ as a sum of \un{two} roots, the summand containing $\ap_j$ lies
in $I$, while the other summand is in $\Delta^+\setminus I$.
This means that $\gamma\not\in I^2$,
i.e., $l(\gamma,I)=1$, and 
$k(\gamma, I)\ge 3$. Thus, condition (ii) in Theorem~\ref{pereform}
is not satisfied, and therefore $I$ is not minimax.
\end{proof*}%
\begin{s}{Theorem}   \label{th-nm-sl}
The number of ideals in $\Delta^+({\frak sl}_{n+1})$ with $k$
non-meeting
generators is equal to $\genfrac{(}{)}{0pt}{}{n}{2k}C_k$.
In particular, the total number of ideals with non-meeting generators
is $M_n$.
\end{s}\begin{proof*}
Let $I$ be an ideal with non-meeting generators and
$\Gamma(I)=\{(a_1,b_1),\ldots,(a_k,b_k)\}$.
In order to compute the number of such ideals,
it is convenient to replace $b_j$ with $\tilde b_j=b_j-1$.
In view of Definition~\ref{not-meet}, the numbers $\{a_i, \tilde b_j\mid
1\le i,j\le n\}$, which belong to $\{1,2,\ldots,,n\}$, 
are pairwise different. They also satisfy the conditions
\begin{equation}  \label{condi}
a_i<\tilde b_i,\quad  a_1<\ldots < a_k,\quad 
\tilde b_1<\ldots < \tilde b_k \ .
\end{equation}
To obtain such a collection of numbers, one should first choose 
arbitrarily $2k$ numbers from $\{1,2,\ldots,n\}$. 
Next, given $2k$ numbers, one should choose
$k$ numbers among them and call them $a_1,\ldots,a_k$
(in the increasing order!). Then the remaining $k$ numbers form the
sequence of $\tilde b_j$'s. But the choice of the $a_i$'s cannot be
arbitrary, for Eq.~\re{condi} must be satisfied.
It is easily seen that the number of admissible choices is
$C_k$. Indeed, given an ordered sequence of $2k$ elements
$v_1\, v_2\, \ldots v_{2k}$ from $\{1,2,\ldots,n\}$,
we assign to $v_l$  the value $+1$,
if $v_l=a_i$; and $-1$, if $v_l=\tilde b_j$. Then Eq.~\re{condi}
is satisfied if and only if all partial sums $\sum_{i\le m}v_i$
are nonnegative. But the number of such sequences
$v_1\, v_2\, \ldots v_{2k}$ is the $k$-th Catalan number, see
\cite[Ex.\,6.19(r)]{rstan}.
\\
Hence, there are $\sum_{k\ge 0}\genfrac{(}{)}{0pt}{}{n}{2k}C_k$
ideals with non-meeting generators, and it is 
the same number that occurs in the proof of Theorem~\ref{sl}.
\end{proof*}%
\begin{s}{Corollary}  \label{nm-mm-sl}
For ${\frak sl}_{n+1}$, the minimax ideals are precisely the ideals
with non-meeting generators.
\end{s}%

\un{$\g={\frak sp}_{2n}$}. \\
We use a standard matrix model of ${\frak sp}_{2n}$ corresponding
to a Witt basis for alternating bilinear form. 
For this basis of ${\Bbb C}^{2n}$, the algebra ${\frak sp}_{2n}$ has the 
following block form:
\[
{\frak sp}_{2n}=\{\left(\begin{array}{cr} A & B \\ C & D \end{array}
\right)\mid B={\widehat B},\ C={\widehat C},\ D=-{\widehat A} \} \ ,
\]
where $A,B,C,D$ are $n\times n$ matrices and $A\mapsto {\widehat A}$ is the 
transpose relative to the antidiagonal.
If $\ov{\be}$ is the standard Borel subalgebra of ${\frak sl}_{2n}$, then
$\be:=\ov{\be}\cap {\frak sp}_{2n}$ is a Borel subalgebra of ${\frak sp}_{2n}$.
(See also \cite[5.1]{duality}.) 
We identify the positive roots of ${\frak sp}_{2n}$ with 
the set $\{(i,j)\mid i<j,\ i+j\le 2n+1\}$.
Here the simple roots are $\ap_i=(i,i+1)$, $1\le i\le n$, and therefore:
\\[.5ex]\indent
\hphantom{$\bullet$ --} $(i,j)=\left\{\begin{array}{cc}
\ap_i+\ldots+\ap_{j-1}, & \text{if }\ j\le n+1 \ ,\\
\ap_i+\ldots +\ap_{2n-j}+2(\ap_{2n-j+1}+\ldots +\ap_{n-1})+\ap_n,
& \text{if }\ j > n+1 \ .\end{array}\right.$
\\[.7ex]
The ideals for ${\frak sp}_{2n}$ can be identified with the ideals 
for ${\frak sl}_{2n}$ that are 
symmetric with respect to the antidiagonal (=\,{\it self-conjugate}). 
In other words, there is a natural bijection between the ideals in
$\Delta^+({\frak sp}_{2n})$ and the self-conjugate ideals in
$\Delta^+({\frak sl}_{2n})$. More precisely, 
suppose $\bar I\in \AD({\frak sl}_{2n})$ and
$\Gamma(\bar I)=\{(i_1,j_1),\ldots,(i_k,j_k)\}$
with $i_1<\ldots < i_k$, where we use our convention on the roots of
${\frak sl}_{2n}$. Then $\bar I$ is
self-conjugate if and only if $i_m+j_{k+1-m}=2n+1$ for all $m$.
The corresponding ideal $I\in\AD({\frak sp}_{2n})$ 
has the generators
$\Gamma(I)=\{(i_m,j_m)\mid m\le [(k+1)/2]\}$. 
We shall say that $\bar I\in\AD({\frak sl}_{2n})$ is the {\it symmetrization\/}
of $I\in \AD({\frak sp}_{2n})$.

%
\begin{s}{Proposition}    \label{nm-sp}
If $I\in\AD({\frak sp}_{2n})$ is a minimax ideal, then 
the symmetrization $\bar I$ has non-meeting 
generators.
\end{s}\begin{proof*}
Suppose the symmetrization
$\bar I$ has some meeting generators. Then arguing as in Proposition~\ref{nm-sl},
we find a root $\bar \gamma\in\bar I$
such that $l(\bar \gamma,\bar I)=1$ and $k(\bar \gamma,\bar I)\ge 3$.
To any root $\bar\gamma\in 
\Delta^+({\frak sl}_{2n})$, 
one naturally associate the root $\gamma\in \Delta^+({\frak sp}_{2n})$. 
With our identification for 
both root systems, this can be formalized as follows.
If $\bar\gamma=(i,j)$ and $i+j\le 2n+1$, then set $\gamma=\bar\gamma$.
If $i+j\ge 2n+2$, then set $\gamma=(2n{+}1{-}j, 2n{+}1{-}i)$. This yields a 
surjective mapping $\Delta^+({\frak sl}_{2n})\to\Delta^+({\frak sp}_{2n})$
and, in particular, $\bar I\to I$ for any $I\in\AD({\frak sp}_{2n})$.
The last and easy observation is that $(\bar I)^l=\ov{I^l}$ for any $I\in 
\AD({\frak sp}_{2n})$.
So that $l(\gamma,I)=1$ and $k(\gamma,I)\ge 3$ as well. Thus,
$I$ is not minimax.
\end{proof*}%
\begin{s}{Theorem}   \label{th-nm-sp}
The number of ideals in $\Delta^+({\frak sp}_{2n})$ with $q$ generators, whose
symmetrization has non-meeting generators,
is equal to 
$\genfrac{(}{)}{0pt}{}{2q-1}{q-1}\genfrac{(}{)}{0pt}{}{n-1}{2q-1}
+\genfrac{(}{)}{0pt}{}{2q}{q}\genfrac{(}{)}{0pt}{}{n-1}{2q}$.
In particular, the total number of such ideals is
$\sum_{q\ge 0}
\genfrac{(}{)}{0pt}{}{q}{[q/2]}\genfrac{(}{)}{0pt}{}{n-1}{q}=\text{dir}_n$.
\end{s}\begin{proof*}
Let $\bar I\in \AD({\frak sl}_{2n})$ be a self-conjugate ideal and  
$\Gamma(\bar I)=\{((i_1,j_1),\ldots,(i_k,j_k)\}$ the sequence of its
generators.
Then 
\[  
i_m<j_m,\quad 1\le i_1< i_2<\ldots < i_k, \quad j_1< j_2<\ldots < j_k\le 2n\ ,
\]  
and the symmetry condition $i_m+j_{k+1-m}=2n+1$ is satisfied for any $m$.
If the generators do not meet, then all 
the numbers $\{i_l, {\tilde j}_m=j_m-1\mid 1\le l,m \le k\}$ are different. 
They form a set $E\subset \{1,2,\ldots, 2n{-}1\}$ consisting of $2k$ 
elements. Because of the symmetry, $E$ is completely determined
by the first $k$ elements, which lie in
$\{1,2,\ldots, n\}$. Moreover, the symmetry and ``non-meeting'' condition 
readily imply
that $n \not\in E$. Thus, $\frac{1}{2}E:=E\cap\{1,\dots,n\}$
actually belongs
to $\{1,2,\ldots, n-1\}$ and $\#(\frac{1}{2}E)=k$. 
Notice that $E$ (and hence $\frac{1}{2}E$) arises as  
a disjoint union of its $i$-part and $j$-part.
So, the problem is to count the admissible partitions in two parts
of all $k$-element
subsets of $\{1,\dots,n{-}1\}$.
To this end, one should first choose arbitrarily $k$ numbers
from $\{1,2,\dots,n{-}1\}$, and then to choose a partition of this set into
$i$- and $j$-parts. In order to compute the number of admissible
partitions, we restate the problem,
as in the proof of Theorem~\ref{th-nm-sl}, 
in terms of sequences of $+1$ and $-1$: Given a $k$-element
subset $v_1<v_2<\ldots <v_k$ of $\{1,\dots,n{-}1\}$, 
we assign the value $+1$ to all elements lying in the $i$-part,
and  $-1$ to all elements lying in the $j$-part. 
It is easily seen that such a sequence corresponds to $\frac{1}{2}E$ 
for a suitable subset $E$ if and only if all partial sums $\sum_{s\le m}v_s$
are non-negative. It is not hard to prove (e.g. using a lattice path 
interpretation
and the reflection principle) that the
number of such sequences is equal to
$\genfrac{(}{)}{0pt}{}{k}{[k/2]}$.
It remains to remember that the ideal $I\in \AD({\frak sp}_{2n})$,
corresponding to $\bar I$,
has $[(k{+}1)/2]$ generators, so that the ideals
with $q$ generators arise if $k=2q{-}1,\,2q$.
\end{proof*}%
\begin{s}{Corollary}  \label{nm-mm-sp}
For $I\in\AD({\frak sp}_{2n})$, the following conditions are equivalent:
\begin{itemize}
\item[\sf (i)] \ $I\in\AD_{mm}({\frak sp}_{2n})$;
\item[\sf (ii)] \ $\bar I\in\AD_{mm}({\frak sl}_{2n})$;
\item[\sf (iii)] \ $\bar I$ has non-meeting generators.
\end{itemize}
\end{s}%
\vskip-1.2ex
\begin{rem}{Example}
By Theorem~\ref{th-nm-sp} and Corollary~\ref{nm-mm-sp}, the number of
minimax ideals with one generator is equal to $(n-1)^2$. It is easy to
verify that the set of positive roots occurring in this way
is $\Delta^+({\frak sp}_{2n})\setminus 
(\Pi\cup\{\ap_i+\ldots +\ap_n\mid i=1,2,\dots,n-1\})$.
\end{rem}%
In \cite{duality}, we considered the statistic on $\AD$ which assigns to
an ideal the number of its generators.
The corresponding generating functions
(polynomials) turn out to be always palindromic. 
It is also makes sense to compute the respective generating functions for
various classes of ideals.
Theorems~\ref{th-nm-sl} and \ref{th-nm-sp}
give us essentially these generating functions for $\AD_{mm}({\frak sl}_{n+1})$ and
$\AD_{mm}({\frak sp}_{2n})$:
\begin{gather*}
{\mathcal F}_{mm}({\frak sl}_{n+1};t)=\sum_{k\ge 0}\genfrac{(}{)}{0pt}{}{n}{2k}C_k t^k
\ , \\
{\mathcal F}_{mm}({\frak sp}_{2n};t)=\sum_{k\ge 0}
(\genfrac{(}{)}{0pt}{}{2k-1}{k-1}\genfrac{(}{)}{0pt}{}{n-1}{2k-1}
+\genfrac{(}{)}{0pt}{}{2k}{k}\genfrac{(}{)}{0pt}{}{n-1}{2k}) t^k \ .
\end{gather*}
It would be interesting to compute the polynomials 
${\mathcal
F}_{mm}(\g;t)$ for all simple Lie algebras. 
It is likely that these polynomials coincide for
${\frak sp}_{2n}$ and ${\frak so}_{2n+1}$, but I have 
no suggestion for ${\frak so}_{2n}$.


\section{Concluding remarks}  
\label{concl}
\setcounter{equation}{0}

\noindent
Here we state several questions/problems related to minimax elements.

1.  Is there a uniform expression for $\#(\HW_{mm})$ for all simple
Lie algebras?

2. Consider the set 
$\displaystyle \bigcup_{w\in \HW_{mm}} w^{-1}\asst\ov{\gA}\subset V$.
It is just the union of the closures of all dominant regions
consisting of a single alcove.
Is there a reasonable description of this set? Note that it is not convex 
in general.

3. Describe combinatorial properties of the polytope $D_{mm}$ defined in 
Proposition~\ref{opis-mm}. Compute the Ehrhart quasi-polynomial for
$D_{mm}$.

4. Find a uniform proof of Theorem~\ref{P:Q}. It is worth noting that the
similar statement can be proved a priori for the simplices $D_{min}$ and 
$D_{max}$. Unfortunately, this does not immediately imply the validity of 
Theorem~\ref{P:Q} for $D_{mm}=D_{min}\cap D_{max}$.

\end{document}